\documentclass[12pt]{article}
\usepackage{amsmath,amssymb}
\usepackage{verbatim}
\usepackage{dsfont,bm}
\usepackage{color}
\usepackage{graphicx}
\usepackage{amsthm}
\usepackage{enumerate}

\usepackage{xcolor}
\definecolor{darkgreen}{rgb}{0,0.75,0}
\definecolor{darkred}{rgb}{0.75,0,0}
\definecolor{darkmagenta}{rgb}{0.5,0,0.5}
\usepackage[colorlinks,citecolor=darkgreen,linkcolor=darkred,urlcolor=darkmagenta]{hyperref}

\newtheorem{theorem}{Theorem}[section]

\newtheorem{lemma}[theorem]{Lemma}

\newtheorem{proposition}[theorem]{Proposition}

\newtheorem{definition}[theorem]{Definition}

\newtheorem{assumption}[theorem]{Assumption}
\newtheorem{remark}[theorem]{Remark}
\newtheorem{example}[theorem]{Example}

\numberwithin{equation}{section}

\def\be{\begin{equation}}
\def\ee{\end{equation}}
\def\bes{\begin{equation*}}
\def\ees{\end{equation*}}

\newcommand{\Cpc}[0]{\operatorname{Cap}}
\newcommand{\mr}[1]{{\tt \href{http://www.ams.org/mathscinet-getitem?mr=#1}{MR#1}}}
\newcommand{\arxiv}[1]{{\tt \href{http://arxiv.org/abs/#1}{arXiv:#1}}}
\newcommand{\set}[1]{\left\{ #1 \right\}}

\newcommand{\Sett}[2]{\left\{ #1  : \, #2 \right\}}

\newcommand{\abs}[1]{{\left\vert\kern-0.25ex #1
    \kern-0.25ex\right\vert}}
\newcommand\norm[1]{\left\lVert#1\right\rVert} 
\newcommand{\one}{\mathds{1}} 
\newcommand{\loc}[0]{\operatorname{loc}}
\newcommand{\diag}[0]{\operatorname{diag}}
\DeclareMathOperator*{\esssup}{ess\,sup}
\DeclareMathOperator*{\essinf}{ess\,inf}

\newcommand{\compl}{c}

  \def\sC {{\mathcal C}}
\def\sD {{\mathcal D}} \def\sE {{\mathcal E}} \def\sF {{\mathcal F}}
  
  \def\sL {{\mathcal L}}
\def\sM {{\mathcal M}}

  \def\sX {{\mathcal X}}

 \def\bE {{\mathbb E}} 
\def\bG {{\mathbb G}}  
  
 \def\bN {{\mathbb N}} 
\def\bP {{\mathbb P}}  \def\bR {{\mathbb R}}
  
\def\bV {{\mathbb V}}

\def\sms{\smallskip}

\def\sm{\smallskip\noindent}

\def\ignore#1{}

\def\ol{\overline}           

\def\eps{\varepsilon}
 
\def\vp{\varphi}
\def\Gam{\Gamma} \def\gam{\gamma}

\def\to {\rightarrow}

\def\pd {\partial}
\def\q{\quad} 
\def\dint{\int\kern-.6em\int}

\def\diam{\mathop{{\rm diam }}}

\def\supp{\mathop{{\rm supp}}}

\def \half {{\tfrac12}}
\def\fract{\tfrac}

\def\wt{\widetilde}

\def\be{\begin{equation}}
\def\ee{\end{equation}}
\def\bes{\begin{equation*}}
\def\ees{\end{equation*}}
\def\ba{\begin{align}}
\def\ea{\end{align}}
\def\xxea{\end{align}}
\def\bas{\begin{align*}}
\def\eas{\end{align*}}

\def\nn{\nonumber} 

\def\proof{{\smallskip\noindent {\em Proof. }}}
\def\qed{{\hfill $\square$ \bigskip}}

\begin{document}

\font\titlefont=cmbx14 scaled\magstep1
\title{\titlefont Boundary Harnack principle and elliptic Harnack inequality}

\author{
M. T. Barlow\footnote{Research partially supported by NSERC (Canada)}, \,
 M. Murugan\footnote{Research partially supported by NSERC (Canada) 
 and the Pacific Institute for the Mathematical Sciences}
}

\maketitle


\begin{abstract}
	We prove a scale-invariant boundary Harnack principle for inner uniform domains over a large family of  Dirichlet spaces. A novel feature of our work is that 
we do not assume volume doubling property for the symmetric measure.
	
\end{abstract}

\section{Introduction}

Let $(\sX,d)$ be a metric space, and assume that associated with this space
is a structure which gives a family of harmonic functions on domains $D \subset \sX$.
(For example, $\bR^d$ with the usual definition of harmonic functions.)
The {\em elliptic Harnack inequality} (EHI) holds if there exists a constant
$C_H$ such that, whenever $h$ is non-negative and harmonic in a ball $B(x,r)$,
then, writing $\half B = B(x,r/2)$,
\be
  \esssup_{\half B} h \le C_H \essinf_{\half B} h. 
\ee
Thus the EHI controls harmonic functions in a domain $D$ away from the boundary $\pd D$.
On the other hand, the {\em boundary Harnack principle} (BHP) controls the ratio of 
two positive harmonic functions near the boundary of a domain. 
The BHP given in \cite{Anc} states that if $D\subset  \bR^d $ is a Lipschitz domain, 
$\xi \in \pd D$, $r>0$ is small enough, then for any pair $u,v$ of non-negative harmonic 
functions in $D$ which vanish on $\pd D \cap B(\xi,2r)$, 
\be
  \frac {u(x)}{v(x)} \le C  \frac {u(y)}{v(y)} \q \text{ for } x,y \in D \cap B(\xi,r). 
\ee  
The BHP is a key component in understanding the behaviour
of harmonic functions near the boundary. It
will in general lead to a characterisation of the Martin boundary, and there
is a close connection between BHP and a Carleson estimate -- see \cite{ALM, Aik08}.
(See also \cite{Aik08} for a discussion of different kinds of BHP.)

The results in \cite{Anc} have been extended in several ways.
The first direction has been to weaken the smoothness hypotheses on the domain $D$;
for example \cite{Aik01} proves a BHP for uniform domains in Euclidean space. 
A second direction is to consider functions which are harmonic
with respect to more general operators. The standard Laplacian is the
(infinitesimal) generator of the semigroup for Brownian motion, and it is
natural to ask about the BHP for more general Markov processes,
with values in a metric space $(\sX,d)$.
In \cite{GyS} the authors prove a BHP for inner uniform domains in a measure
metric space $(\sX,d,m)$ with a Dirichlet form which satisfies the
standard parabolic Harnack inequality (PHI). 
These results are extended in \cite{Lie} to spaces which satisfy a 
parabolic Harnack inequality with anomalous space-time scaling.
In most cases the BHP has been proved for Markov processes which are
symmetric, but see \cite{LS} for the BHP for some more general processes.
All the papers cited above study  the harmonic functions associated with
continuous Markov processes: see \cite{Bog,BKK} for a BHP for a
class of jump processes. 

The starting point for this paper is the observation that the BHP
is a purely elliptic result, and one might expect that the proof would only
use elliptic data. However, the generalizations of the BHP 
beyond the Euclidean case in \cite{GyS, LS, Lie} all 
use parabolic data, or more precisely bounds on the heat kernel of the process. 



The main result of this paper is as follows. See Sections \ref{sec:iud},
\ref{sec:Dspace} and \ref{sec:harm} for unexplained definitions and notation.

\begin{theorem}\label{T:BHP}
Let $(\sX,d)$ be a complete, separable, locally compact, 
length space, and let $\mu$ be a non atomic Radon measure
on $(\sX,d)$ with full support. 
Let $(\sE,\sF)$  be a regular strongly local Dirichlet form on $L^2(\sX,\mu)$. 
Assume that $(\sX,d,\mu,\sE,\sF)$ satisfies the 
elliptic Harnack inequality, and has Green functions which satisfy
the regularity hypothesis Assumption \ref{a-green}. 
Let $U \subsetneq \sX$ be an inner uniform domain.
Then there exist $A_0,C_1 \in (1,\infty)$, $R(U) \in (0,\infty]$ such that for all  
$\xi \in \partial_{\widetilde{U}}U$, for all $0<r<R(U)$ and any two non-negative functions 
$u,v$ that are harmonic on $B_U(\xi,A_0r)$ with Dirichlet boundary conditions along 
$\partial_{\widetilde{U}}U \cap B_{\widetilde{U}}(\xi,2A_0r)$, we have
\[
 \esssup_{x \in B_U(\xi,r)}	\frac{u(x)}{v(x)} 
 \le C_1  \essinf_{x \in B_U(\xi,r)} \frac{u(x)}{v(x)} .
\]
The constant $R(U)$ depends only on the inner uniformity constants of $U$ and 
$\operatorname{diameter}(U)$, and can be chosen to be $+\infty$ if $U$ is unbounded.
\end{theorem}

\begin{remark}
{\rm (1) The constant $A_0$ above depends only on the inner uniformity 
constants for the domain $U$, and 
$C_1$ depends only on these constants and the constants  in the EHI. \\
(2) Since the EHI is weaker than the PHI, our result extends the BHP
to a wider class of spaces; also our approach has the advantage that
we can dispense with heat kernel bounds. 
Our main result provides new examples of differential operators 
that satisfy the BHP even in $\mathbb{R}^n$ -- see \cite[ (2.1) and Example 6.14]{GS}. \\
(3) By the standard oscillation lemma  (see \cite[Lemma 5.2]{GT12}), any locally bounded harmonic function admits a continuous version. The elliptic Harnack inequality implies that any non-negative harmonic function is locally bounded. Therefore, under our assumptions, every non-negative harmonic function admits a continuous version. \\ 
(4) Let $\mu'$ be a measure which is mutually absolutely continuous
with respect to the measure $\mu$ in the Theorem above, and suppose
that $d\mu'/d\mu$ is bounded away from 0 and infinity on compact subsets of $\sX$. Then 
(see Remark \ref{R:com}) this change of measure does not change the family of harmonic
functions, or the Green functions, and the
hypotheses of Theorem \ref{T:BHP} hold for  $(\sX,d,\mu',\sE,\sF)$. 
On the other hand,  heat kernel bounds and parabolic Harnack inequality are not in general preserved by
such a change of measure because $d\mu'/d\mu$ need not be bounded away from 0 or infinity on $\sX$. 
}\end{remark}

%

The contents of this paper are as follows. In Section \ref{sec:iud} we give the 
definition and basic properties of inner uniform domains in length spaces. 
Section \ref{sec:Dspace} reviews the properties of Dirichlet forms and
the associated Hunt processes. In Section \ref{sec:harm} we give the
definition of harmonic function in our context, state 
Assumption \ref{a-green}, and give some consequences. In particular,
we prove the essential technical result that Green
functions are locally in the domain of the Dirichlet form -- see Lemma \ref{L:gdom}. 
The key comparisons of Green functions, which follow from the EHI,
and were proved in \cite{BM}, are given in Proposition \ref{p-consq}.
In the second part of this section we give some sufficient 
conditions for Assumption \ref{a-green} to hold, in terms of local ultracontracivity.
We conclude Section \ref{sec:harm} with two examples: weighted manifolds
and cable systems of graphs.

After these rather lengthy preliminaries, Section \ref{sec:bhp} gives
the proof of Theorem \ref{T:BHP}. 
We follow Aikawa's approach in \cite{Aik01},
which proved the BHP for uniform domains in $\bR^n$. 
This method has been
adapted to more general settings \cite{ALM, GyS, LS, Lie}. 
The papers \cite{GyS, LS, Lie} all consider domains in more general metric spaces,
and use heat kernel estimates to 
obtain two sided estimates for the Green function in a domain; these
estimates are then used in the proof of the BHP.
For example \cite[Lemma 4.5]{Lie} gives upper and lower bounds 
on $g_D(x,y)$ when $D$ is a domain of diameter $R$, and the points
$x,y$ are separated from $\pd D$ and each other by a distance greater than
$\delta R$. These bounds are of the form
$\Psi(R)/\mu(B(x,R))$; here $\Psi: [0,\infty) \to [0,\infty)$ is a global space time
scaling function. (See \cite{Lie} for the precise statement.)
In our argument we use instead the comparison of Green functions
given by Proposition \ref{p-consq}.

We use $c, c', C, C'$ for strictly positive constants, which may change value 
from line to line. Constants with numerical subscripts will keep the same value in each
argument, while those with letter subscripts will be regarded as constant
throughout the paper.
The notation $C_0 = C_0(a)$ means that the constant
$C_0$ depends only on the constant $a$.

\section{Inner uniform domains} \label{sec:iud}

In this section we introduce the geometric assumptions on the underlying 
metric space, and the corresponding domains.

\begin{definition}[Length space]
{\rm 
Let $(\sX,d)$ be a  metric space.
The length  $L(\gamma) \in [0,\infty]$ of a continuous curve $\gamma:[0,1] \to \sX$ is given by
\[
 L(\gamma)= \sup  \sum_{i} d( \gamma(t_{i-1}),\gamma(t_i)),
\]
where the supremum is taken over all partitions $0=t_0 < t_1 < \ldots < t_k = 1$ of $[0,1]$.
It is clear that $L(\gamma) \ge d(\gamma(0),\gamma(1))$. 
A metric space is a \emph{length space} if $d(x,y)$ is equal to
the infimum of the lengths of continuous curves joining $x$ and $y$.
}\end{definition}
 
 For the rest of this paper we will assume that $(\sX,d)$ is
a complete, separable, locally compact, length space.
We write $\ol A$ and $\pd A$ for the closure and boundary respectively of a subset
$A$ in $\sX$.
By the Hopf--Rinow--Cohn-Vossen theorem (cf. \cite[Theorem  2.5.28]{BBI}) 
every closed metric ball in $(\sX,d)$ is compact. 
It also follows that there exists a geodesic path $\gam(x,y)$ (not necessarily unique) between any two
points $x,y \in \sX$.
We write $B(x,r)= \Sett{y \in \sX} {d(x,y) < r}$ for open balls in $(\sX,d)$.

Next, we introduce the intrinsic distance $d_U$ induced by an open set $U \subset \sX$.

\begin{definition}[Intrinsic distance]
{\rm 
Let $U \subset \sX$ be a connected open subset. We define the \emph{intrinsic distance} $d_U$ by
\[
d_U(x,y) = \inf \Sett{L(\gamma)}{\gamma:[0,1] \to U \mbox{ continuous, } \gamma(0)= x, \gamma(1)=y}.
\]
}\end{definition}

It is well-known that $(U,d_U)$ is a length space (cf. \cite[Exercise 2.4.15]{BBI}). 
We now consider its completion.

\begin{definition}[Balls in intrinsic metric]
{\rm Let $U \subset \sX$ be connected and open.
 Let $\widetilde{U}$ denote the completion of $(U,d_U)$, equipped 
with the natural extension of  $d_U$ to $\widetilde{U} \times \widetilde{U}$.
For $x \in \wt U$ we define 
\[
B_{\widetilde U}(x,r) = \Sett{y \in \widetilde U}{d_U(x,y) < r}.
\]
Set 
\[
B_U(x,r)= U \cap B_{\widetilde U}(x,r).
\]
}\end{definition}
If $x \in U$, then $B_U(x,r)$ simply corresponds to the open ball in $(U,d_U)$. However, 
the definition of $B_U(x,r)$ also makes sense for $x \in \widetilde U \setminus U$.

\begin{definition}[Boundary and distance to the boundary]
{\rm We denote the boundary of $U$ with respect to the inner metric by
\[
\partial_{\widetilde U} U = \widetilde U \setminus U,
\]
and the distance to the boundary by
\[
\delta_U(x) = \inf_{y \in \partial_{\widetilde U} U} d_U(x,y) = \inf_{ y \in \sX \setminus U} d(x,y).
\]
For any open set $V \subset U$, let $\overline V^{d_U}$ denote the completion of $V$ 
with respect to the metric $d_U$. We denote the boundary of $V$ with respect to 
$\widetilde U$ by
\[
\partial_{\widetilde U} V = \overline V^{d_U} \setminus V.
\]
}\end{definition}


\begin{definition}[Inner uniform domain]
\label{d:iud}
{\rm Let $U$ be a connected, open subset of a length space $(\sX,d)$.
Let $\gamma:[0,1] \to U$ be a rectifiable, continuous curve in $U$. Let $c_U,C_U \in (0,\infty)$. 
We say $\gamma$ is a {\em $(c_U,C_U)$-inner uniform curve} if
\[
L(\gamma) \le C_U d_U(\gamma(0),\gamma(1)),
\]
and
\[
\delta_U(\gamma(t)) \ge c_U \min \left( d_U(\gamma(0), \gamma(t)), d_U(\gamma(1),\gamma(t)) \right)
\mbox{ for all $t \in [0,1]$.}
\]
The domain $U$ is called a \emph{$(c_U,C_U)$-inner uniform domain} if any two points in $U$ 
can be joined by a $(c_U,C_U)$-inner uniform curve.
%
}\end{definition}

%
%
%

The following lemma extends the existence of inner uniform curves between any 
two points in $U$ in Definition \ref{d:iud} to  the existence of inner uniform curves 
between  any two points in $\wt{U}$.

\begin{lemma} \label{l:iuc}
Let $(\sX,d)$ be a complete, locally compact, separable, length space. 
Let $U$ be a  $(c_U,C_U)$-inner uniform domain and let $\wt{U}$ denote the 
completion of $U$ with respect to the inner metric $d_U$. Then for any two 
points $x,y$ in $(\wt{U},d_U)$, there exists a $(c_U,C_U)$-uniform curve in the $d_U$ metric.
\end{lemma}

\proof Let $x,y \in \wt{U}$.
There exist sequences $(x_n),(y_n)$ in $U$ such that $x_n \to y, y_n \to y$ 
in the $d_U$ metric as $n \to \infty$.
Let $\gamma_n:[0,1] \to U, n \in \bN$ be a $(c_U,C_U)$-inner uniform curve 
in $U$ from $x_n$ to $y_n$ with constant speed parametrization. 
By \cite[Theorem  2.5.28]{BBI}, 
the curves $\gam_n$ can be viewed as being in the compact space  
$\overline{B_U(x,2C_U d_U(x,y))}^{d_U}$ for all large enough $n$. 
By a version of Arzela-Ascoli theorem the desired inner uniform curve 
$\gamma$ from $x$ to $y$ can be constructed as a sub-sequential limit of the curves
$(\gam_n)$ -- see \cite[Theorem 2.5.14]{BBI}.
\qed

The following geometric property of a  metric space $(\sX,d)$ will play an important
role in the paper.

\begin{definition}[Metric doubling property]
{\rm We say that a metric space $(\sX,d)$ satisfies the \emph{metric doubling property} 
if there exists $C_M>0$ such that  any ball $B(x,r)$ can be covered by at most
$C_M$ balls of radius $r/2$. 
}\end{definition}

%


%

Let $\overline U \subset \sX$ denote the closure of $U$ in $(\sX,d)$.
Let $p: (\widetilde U, d_U) \to ( \overline{U}, d)$ denote the natural projection map, that is $p$ is the unique continuous map such that $p$ restricted to $U$ is the identity map on $U$. 
The following lemma  allows us to compare balls with respect to the
$d$ and $d_U$ metrics.

\begin{lemma}\label{l:ballcomp}
Let $(\sX,d)$ be a complete, length space satisfying the metric doubling property. Let $U \subset \sX$ be a connected, open, $(c_U,C_U)$-inner uniform domain.
Then there exists $ \widetilde{C_U} >1$ such that 
for all balls $B(p(x),r/ \widetilde{C_U}) $ with $x \in \widetilde U$ and $ r >0$, we have
\[
B_{\widetilde U}(x,r/ \widetilde{C_U}) \subset D' \subset B_{\widetilde U}(x,r),
\]
where $D'$ is the  connected component of $p^{-1}( B(p(x),r/\widetilde{C_U})  \cap \overline U)$ containing $x$.
\end{lemma}

\proof
See \cite[Lemma 3.7]{LS} where this is proved under the hypothesis
of volume doubling, and note that the argument only uses metric doubling.
(Alternatively, a doubling measure exists by \cite[Theorem 1]{LuS}, and one can then use
\cite{LS}). \qed

The following lemma shows that every point in an inner uniform domain is close to a point that is sufficiently far away from the boundary.

\begin{lemma}(\cite[Lemma 3.20]{GyS})\label{l:faraway}
Let $U$ be a $(c_U,C_U)$-inner uniform domain in a length space $(\sX,d)$. 
For every inner ball $B=B_{\widetilde U}(x,r)$ with the property that
$B\neq B_{\widetilde U}(x,2r)$ there exists a point $x_r \in B$ with
\[
d_U(x,x_r)=r/4 \hspace{2mm} \mbox{ and }\hspace{2mm} \delta_U(x_r) \ge \frac{c_U r }{4}.
\]
\end{lemma}


\begin{lemma} \label{L:goodcurve}
Let $U$ be a $(c_U,C_U)$-inner uniform domain in a length space $(\sX,d)$. 
If $x, y \in U$, then there exists a $(c_U,C_U)$-inner uniform curve
$\gam$ connecting $x$ and $y$ with 
$\delta_U(z) \ge \half c_U\left(  \delta_U(x) \wedge  \delta_U(y) \right)$ for all $z \in \gam$.
\end{lemma}

\proof Write $t=\delta_U(x) \wedge  \delta_U(y)$. Let $\gam$ be an inner uniform
curve from $x$ to $y$ and let $z \in \gam$.
If $d_U(x,z) \le \half t$, then $\delta_U(z) \ge \delta_U(x) - d_U(x,z)\ge \half t$,
and the same bound holds if $d_U(y,z) \le \half t$. Finally 
if $d_U(x,z) \wedge d_U(y,z) \ge \half t$, then $\delta_U(z) \ge \half c_U t$. \qed

\section{Dirichlet spaces and Hunt processes} \label{sec:Dspace}

Let $(\sX,d)$ be a locally compact, separable, metric space and let $\mu$ be a Radon measure with full support.
Let $(\sE,\sF)$ be a regular, strongly local Dirichlet form on $L^2(\sX,\mu)$
-- see \cite{FOT}. Recall that a Dirichlet
form $(\sE,\sF)$ is \emph{strongly local} if $\sE(f,g) =0$ for any $f,g \in \sF$ with 
compact supports, such that $f$ is constant in a neighbourhood of $\supp(g)$.
We call $(\sX, d,\mu,\sE,\sF)$ 
a {\em metric measure Dirichlet} space, or {\em MMD space} for short.

Let $\sL$ be the generator of $(\sE,\sF)$ in $L^2(\sX,\mu)$; that is $\sL$ is a self-adjoint and non-positive-definite operator in $L^2(\sX,\mu)$ with domain $\sD(\sL)$ that is dense in $\sF$ such that
\[
\sE(f,g) = - \langle \sL f, g\rangle,
\]
for all $f \in \sD(\sL)$ and for all $g \in \sF$;
here $\langle \cdot, \cdot\rangle$, is the inner product in $L^2(\sX,\mu)$. 
The associated \emph{heat semigroup}
\[
P_t = e ^{t \sL}, t \ge 0,
\]
is a family of contractive, strongly continuous, Markovian, self-adjoint operators in $L^2(\sX,\mu)$. 
We set
\be
 \sE_1(f,g) = \sE(f,g) + \langle f,g\rangle, \q
  ||f||_{\sE_1} = \sE_1(f,f)^{1/2}.
\ee

It is known that corresponding to a regular Dirichlet form, there exists an essentially 
unique Hunt process $X=(X_t, t\ge 0, \bP^x, x\in \sX)$.
The relation between the Dirichlet form $(\sE,\sF)$ on $L^2(\sX,\mu)$ and the associated Hunt 
process is given by the  identity
\[ P_tf(x) = \bE^x f(X_t), \]
for all $f \in L^\infty(\sX,\mu)$, for every $t>0$, and for $\mu$-almost all $x \in \sX$.
Also associated with the Dirichlet form and $f \in \sF$ is the energy measure
$d\Gam(f,f)$. This is defined to be the unique Radon measure 
such that for all $g \in \sF \cap \sC_c(\sX)$, we have 
$$ \int_\sX g \, d\Gam(f,f) = 2 \sE(f,fg) - \sE(f^2,g). $$
We have
$$ \sE(f,f) = \frac{1}{2} \int_{\sX} d\Gam(f,f). $$

\begin{definition}
For an open subset of $U$ of $\sX$, we define the following function spaces 
associated with $(\sE,\sF)$.
\begin{align*}
\sF_{\loc}(U) &= \Sett{u \in L^2_{\loc}(U,\mu)}{ \forall \mbox{ relatively compact open } 
V \subset U, \exists u^\# \in \sF, u = u^\# \big|_{V}  \mbox{ $\mu$-a.e.}}, \\
 \sF(U)   &=  \Sett{u \in \sF_{\loc}(U) }{ \int_U \abs{u}^2 \,d\mu + \int_U d\Gamma(u,u) < \infty }, \\
\sF_{\mbox{\tiny{c}}}(U)  
&=  \Sett{u \in \sF(U)}{\mbox{the essential support of } u \mbox{ is compact in } U }, \\
\sF^0(U)  &= \mbox{ the closure of } \sF_{\mbox{\tiny{c}}}(U) \mbox{ in } \sF \mbox{ in the norm }
||\cdot ||_{\sE_1}.
\end{align*}
\end{definition}

We define capacities for $(\sX,d,\mu,\sE,\sF)$ as follows. 
Let $U$ be an open subset of $\sX$.
By $A \Subset U$, we mean that the closure of $A$ is a compact subset of $U$.
For $A \Subset U$ we set
\be \label{e:capdef}
 \Cpc_U(A) = \inf\{ \sE(f,f): f \in \sF^0(U) \mbox{ and  $f \ge 1$ in a neighbourhood of $A$} \}.
\ee
A statement depending on $x \in A$ is said to hold {\em quasi-everywhere on $A$} 
(abbreviated as q.e.~on $A$), if there exists a set $N \subset A$ of zero capacity 
such that the statement is true for every $x \in A \setminus N$.
It is known that every function $f \in \sF$ admits a quasi continuous version,
which is unique up to a set of zero capacity 
(cf. \cite[Theorem 2.1.3]{FOT}).
Throughout this paper we will assume that every $f \in \sF$ is represented 
by its quasi-continuous version. 

For an open set $U$ an equivalent definition of $\sF^0(U)$ is given by
\begin{equation}\label{e:F0}
\sF^0(U) = \Sett{u \in \sF} { \widetilde{u} = 0 \mbox{ q.e. on $\sX \setminus U$}},
\end{equation}
where $\widetilde{u}$ is a quasi-continuous version of $u$ -- see 
\cite[Theorem 4.4.3(i)]{FOT}. Thus we can 
identify $\sF^0(U)$ as a subset of $L^2(U,\mu)$, where in turn $L^2(U,\mu)$ is 
identified with the subspace $\Sett{u \in L^2(\sX,\mu)}{ u = 0 \mbox{ $\mu$-a.e. on $\sX \setminus U$}}$.

\begin{definition}\label{d-dirdir}
{\rm For an open set $U \subset \sX$, we define the 
{\em part of the Dirichlet form $(\sE,\sF)$ on $U$} by
\[
\mathcal{D}(\sE_U)= \sF^0(U) \mbox{ and } \sE_U(f,g) = \sE(f,g) \mbox{ for $f,g \in \sF^0(U)$}.
\]
By \cite[Theorem 3.3.9]{CF} 
$(\sE_U,\sF^0(U))$ is a regular, strongly local Dirichlet form on $L^2(U,\mu)$.  
We write $(P^U_t, t\ge0)$ for the associated semigroup, and call
$(P^U_t)$ the {\em semigroup of $X$ killed on exiting $U$}.
The Dirichlet form $(\sE_U,\sF^0(U))$ is associated with the process 
$X$ killed upon exiting $U$ -- see \cite[Theorem 3.3.8(ii)]{CF}.
}\end{definition}

For an open set $U$ we need to consider functions that vanish on a portion of the 
boundary of $U$, and we
therefore define the following local spaces associated with $(\sE_U, \sF^0(U))$.

\begin{definition}
{\rm Let $V \subset U$ be open subsets of $\sX$. Set
\begin{align*}
 \sF^0_{\loc}(U,V) = & \{ f \in L^2_{\loc}(V,\mu): \forall \mbox{ open } A \subset V \mbox{ relatively~compact in } \overline{U} \mbox{ with } \\
& \, d_{U}(A,U \setminus V) > 0, \, \exists f^{\sharp} \in \sF^0(U): f^{\sharp} = f \, \mu \mbox{-a.e.~on } A \}.
\end{align*}
}\end{definition}

Note that $\sF^0_{\loc}(U,V) \subset \sF_{\loc}(V)$.
Roughly speaking, a function in $\sF^0_{\loc}(U,V)$ vanishes along the portion of boundary 
given by $\partial_{\widetilde{U}} V \cap \partial_{\widetilde{U}} U$.


\section{Harmonic functions and Green functions}
\label{sec:harm}

\subsection{Harmonic functions}
 
We begin by defining harmonic functions for a strongly 
local, regular Dirichlet form $(\sE,\sF)$ on $L^2(\sX,\mu)$.

\begin{definition}\label{d-harmonic}
{\rm Let $U \subset \sX$ be open. A function $u: U \to \bR$ is \emph{harmonic} on $U$
if  $u \in \sF_{\loc}(U)$ and for any function $\phi \in \sF_{\mbox{\tiny{c}}}(U)$
there exists a function $u^\# \in \sF$ such that $u^\# = u$ in a neighbourhood of the essential support of $\phi$
and
\[
\sE(u^\#, \phi)= 0.
\]
}\end{definition}

\begin{remark}
{\rm 
\begin{enumerate}[(a)]
\item By the locality of $(\sE,\sF)$, $\sE(u^\#, \phi)$ does not depend 
on the choice of $u^\#$ in Definition \ref{d-harmonic}.
\item If $U$ and $V$ are open subsets of $\sX$ with $V \subset U$
and $u$ is harmonic in $U$, then the restriction $u\big|_V$ is harmonic in $V$. 
This follows from the locality of $(\sE,\sF)$.
\item It is known that $u \in L^\infty_{\loc}(U,\mu)$ is harmonic in $U$ if and only if it satisfies the following property: for every relatively compact open subset $V$
of $U$, $t \mapsto \widetilde{u}(X_{t \wedge \tau_V})$ is a uniformly integrable $\bP^x$-martingale 
for q.e.~$x \in V$. (Here $\widetilde{u}$ is a quasi continuous version of $u$ on $V$.)
This equivalence between the weak solution formulation in Definition \ref{d-harmonic} and the probabilistic formulation using martingales is given in \cite[Theorem 2.11]{Che}.
\end{enumerate}
}\end{remark}

\begin{definition}
{\rm Let $V \subset U$ be open. We write $\widetilde{V}^{d_U}$ 
for the closure of $V$ in $(\widetilde{U},d_U)$.
We say that a harmonic function 
$u : V \to \bR$ satisfies \emph{Dirichlet boundary conditions on the boundary} 
$\partial_{\widetilde{U}} U \cap \widetilde{V}^{d_U}$
if $u \in \sF^0_{\loc}(U,V)$.
}\end{definition}

\subsection{Elliptic Harnack inequality}\label{d-ehi}

\begin{definition}[Elliptic Harnack inequality]
{\rm We say that $(\sE,\sF)$ satisfies the \emph{local elliptic Harnack inequality}, denoted 
$\operatorname{EHI}_{\loc}$, if 
there exist constants $C_H <\infty$, $R_0 \in (0, \infty]$ 
and $\delta \in (0,1)$ such that, 
for any ball $B(x,R) \subset \sX$ satisfying $R \in (0,R_0)$, and any non-negative function 
$u \in \sF_{\loc}(B(x,R))$ that is harmonic on $B(x,R)$, we have
\begin{equation*}\label{ehi} 
\esssup_{z \in B(x,\delta R)} u(z) \leq C_H \essinf_{z \in B(x,\delta R)} u(z). 
\tag*{$(\operatorname{EHI})$}
\end{equation*}
We say that  $(\sE,\sF)$ satisfies the \emph{elliptic Harnack inequality},
denoted $(\operatorname{EHI})$, 
if $\operatorname{EHI}_{\loc}$ holds with $R_0=\infty$. 
}\end{definition}


An easy chaining argument along geodesics shows that if the EHI holds 
for some $\delta \in (0,1)$, then it holds for any other $\delta' \in (0,1)$.
Further, if the local EHI holds for some $R_0$, then it holds (with of course a
different constant $C_H$) for any $R'_0 \in (0,\infty)$.

We recall the definition of Harnack chain -- see \cite[Section 3]{JK}. For a ball $B=B(x,r)$, we use the notation $M^{-1} B$ to denote the ball $B(x,M^{-1} r)$.

\begin{definition}[Harnack chain]  
{\rm Let $U \subsetneq \sX$ be a connected open set.
For $x,y \in U$, an {\em $M$-Harnack chain from $x$ to $y$}
in $U$ is a sequence of balls $B_1,B_2,\ldots,B_n$ each contained in $U$ such that $x \in M^{-1}B_1, y \in M^{-1}B_n$, and $M^{-1}B_{i} \cap M^{-1}B_{i+1} \neq \emptyset$, for $i=1,2,\ldots,n-1$. 
The number $n$ of balls in a Harnack chain is called the {\em length} of the Harnack chain.
For a domain $U$ write
$N_U(x,y;M)$ for the length of the shortest $M$-Harnack chain in $U$ from
$x$ to $y$. 
}\end{definition}

\begin{remark} \label{r-chaining}
{\rm Suppose that $(\sE,\sF)$ satisfies the elliptic Harnack inequality
with constants $C_H$ and $\delta$. 
If $u$ is a positive continuous harmonic function on a domain $U$, then
\be \label{e:hchain}
     C_H^{-N_U(x_1, x_2; \delta^{-1})} u(x_1) \le u(x_2) \le C_H^{N_U(x_1, x_2; \delta^{-1})} u(x_1).
\ee
for all $x_1$, $x_2 \in U$.
}\end{remark}

\begin{lemma} 
\label{l:qhmetric}
Let $(\sX,d)$ be a locally compact, separable, length space that satisfies the metric doubling property.
Let $U \subsetneq \sX$ be a $(c_U,C_U)$-inner uniform domain in $(\sX,d)$. 
Then for each $M > 1$ there exists 
$C \in (0,\infty)$,  depending only on $c_U$, $C_U$ and $M$, such that for all $x,y \in U$
\[
C^{-1} \log \left( \frac{d_U(x,y)}{ \min(\delta_U(x),\delta_U(y))} + 1 \right) \le N_U(x,y;M) 
\le C \log \left( \frac{d_U(x,y)}{ \min(\delta_U(x),\delta_U(y))} + 1 \right)  + C.
\]
\end{lemma}

\proof See  \cite[Equation (1.2) and Theorem 1.1]{GO} or \cite[Theorem 3.8 and 3.9]{Aik15}
for a similar statement for the quasi-hyperbolic metric on $U$; the result
then follows by a comparison between the quasi-hyperbolic metric and 
the length of Harnack chains as in \cite[p. 127]{Aik01}.
 \qed

\subsection{Green function}

Let $(\sE, \sF)$ be a regular, strongly local Dirichlet form and 
let $\Omega \subsetneq \sX$ be open. We define 
\[
\lambda_{\min}(\Omega) 
= \inf_{u \in \sF^0(\Omega) \setminus \set{0}} \frac{\sE_\Omega(u,u)}{\norm{u}_2^2}.
\]
Writing $\sL^\Omega$ for the generator of $(\sE_\Omega,\sF^0(\Omega))$, we have
$\lambda_{\min}(\Omega)= \inf \operatorname{spectrum}(-\sL^\Omega)$.

The next Lemma gives the existence and some fundamental properties
of the Green operator on a domain $\Omega \subset \sX$.

\begin{lemma}(\cite[Lemma 5.1]{GH}) \label{l:5.1}
Let $(\sE,\sF)$ be a regular, Dirichlet form in $L^2(\sX,\mu)$ and let $\Omega \subset \sX$ 
be open and satisfy $\lambda_{\min}(\Omega) > 0$. Let $\sL^\Omega$ be the generator of 
$(\sE_\Omega,\sF^0(\Omega))$, and let $G^\Omega = ( - \sL^\Omega)^{-1}$ 
be the inverse of $-\sL^\Omega$ on $L^2(\Omega,\mu)$.
Then the following statements hold:
\begin{enumerate}[(i)]
\item $\norm{G^\Omega} \le \lambda_{\min}(\Omega)^{-1}$, that is, for any $f \in L^2(\Omega,\mu)$,
\[
\norm{G^\Omega f}_{L^2(\Omega)} \le \lambda_{\min}(\Omega)^{-1} \norm{f}_{L^2(\Omega)};
\]
\item for any $f \in L^2(\Omega)$, we have that $G^\Omega f \in \sF^0(\Omega)$, and
\[
\sE_\Omega(G^\Omega f, \phi) = \langle f, \phi \rangle \mbox{ for any $\phi \in \sF^0(\Omega)$};
\]
\item for any $f \in L^2(\Omega)$,
\[
G^\Omega f = \int_0^\infty P_s^\Omega f \, ds;
\]
\item $G^\Omega$ is positivity preserving: $G^\Omega f \ge 0$ if $f \ge 0$.
\end{enumerate}
\end{lemma}

We now state our fundamental assumption on the Green function.

\begin{assumption}\label{a-green}
{\rm Let $(\sX,d)$ be a complete, locally compact, separable, length space 
and let $\mu$ be a non-atomic Radon measure on $(\sX,d)$ with full support.
Let $(\sE,\sF)$ be a strongly local, regular, Dirichlet form on $L^2(\sX,\mu)$. 
Let $\Omega \subset \sX$ be a non-empty bounded open set with 
$\operatorname{diameter}(\Omega,d) \le \operatorname{diameter}(\sX,d)/5$.
Assume that $\lambda_{\min}(\Omega)  > 0$, and
there exists a function $g_\Omega(x,y)$ defined for $(x,y) \in \Omega \times \Omega$ 
with the following properties:
\begin{enumerate}[(i)]
 \item  (Integral kernel) $G^\Omega f(x) = \int_\Omega g_\Omega(x,z) f(x)\, \mu(dz)$ for all $f \in L^2(\Omega)$ and $\mu$-a.e.~$x\in \Omega$;
 \item (Symmetry) $g_\Omega(x,y)= g_\Omega(y,x) \ge 0$ for all $(x,y) \in \Omega \times \Omega \setminus \diag$;
 \item (Continuity) $g_\Omega(x,y)$ is jointly continuous in $(x,y) \in \Omega \times \Omega \setminus \diag$;
 \item (Maximum principles) If $x_0 \in U \Subset \Omega$, then
 \begin{align*}
 \inf_{U \setminus \set{x_0}} g_\Omega(x_0,\cdot) &= \inf_{\partial U} g_\Omega(x_0,\cdot),\\
 \sup_{\Omega \setminus U} g_\Omega(x_0, \cdot) &= \sup_{\partial U}g_\Omega(x_0,\cdot).
 \end{align*}
\end{enumerate}
}\end{assumption}

We now give some consequences of this assumption; 
in the next subsection we will give some sufficient conditions 
for it to hold.

We begin by showing that the Green function $g_\Omega(x,\cdot)$ is harmonic
in $\Omega \setminus \{x\}$ and vanishes along the boundary of $\Omega$. 
Since we are using 
Definition \ref{d-harmonic}, we need first to prove that this function is
locally in the domain of the Dirichlet form, that is that 
$g_\Omega(x,\cdot) \in \sF_{\loc}(\Omega \setminus \set{x})$.
For this it is enough that 
$g_\Omega(x,\cdot) \in \sF^0_{\loc}(\Omega, \Omega \setminus \set{x})$.
This result was shown under more restrictive hypothesis (Gaussian 
or sub-Gaussian heat kernel estimates) in \cite[Lemma 4.7]{GyS} and by similar 
methods in \cite[Lemma 4.3]{Lie}. 
Our proof is based on a different approach (see \cite[Theorem 4.16]{GyS}), using 
Lemma \ref{l:5.1}.

\begin{lemma}\label{L:gdom}
Let $(\sX,d,\mu,\sE,\sF)$, $\Omega$ be as in Assuption \ref{a-green}. 
For  any fixed $x \in \Omega$, the function $y \mapsto g_\Omega(x,y)$ is in
$\sF^0_{\loc}(\Omega, \Omega \setminus \set{x})$, 
and is harmonic in $\Omega \setminus \set{x}$.
\end{lemma}

\proof
Fix $x \in \Omega$.
Let $V \subset \Omega$ be an  open set such that 
$\overline{V} \subset \overline{\Omega} \setminus \set{x}$. 
Let $\Omega_1, \Omega_2$ be precompact open sets such that $\overline{\Omega} \subset \Omega_1 \subset \overline{\Omega_1} \subset \Omega_2$. Let $r>0$ be such that 
$B(x,4r) \subset \Omega \cap V^\compl$.
Let $\phi \in \sF$ be a continuous function such that $0 \le \phi \le 1$ and
\begin{equation*}
\phi = \begin{cases}
1 &\mbox{on $B(x,3r)^\compl \cap \Omega_1$,}\\
0 & \mbox{on $B(x,2r) \cup \left(\overline{\Omega_2}\right)^\compl$}.
\end{cases}
\end{equation*}
Since $\vp \equiv 1$ on $V$, to prove that 
$g_\Omega(x,\cdot) \in \sF^0_{\loc}(\Omega, \Omega \setminus \set{x})$ it is 
sufficient to prove that $\vp g_\Omega(x,\cdot) \in \sF^0(\Omega)$.

For $k \ge 1$ set $B_k = B(x, r/k)$. 
Consider the sequence of functions defined by
\be \label{e:appxG}
h_k(y) = \frac{1}{\mu(B_k)}\int_{B_k} g_\Omega(z,y) \mu(dz), \q 
y \in \Omega, \ k \in \bN.
\ee
By Lemma \ref{l:5.1}(ii) we have $h_k \in \sF^0(\Omega)$ for all $k \ge 1$.

By the maximum principle, we have
\[
M := \sup_{z \in \overline{B(x,r)}, y \in B(x,2r)^\compl} g_\Omega(z,y) = \sup_{z \in \overline{B(x,r)}, y \in \partial B(x,2r)} g_\Omega(z,y) < \infty,
\]
since the image of the compact set $\overline{B(x,r)} \times \partial B(x,2r)$ under the continuous 
map of $g_\Omega$  is bounded. Thus the functions $\phi h_k$ are bounded uniformly
by $M \one_{\Omega \setminus B(x,2r)}$.
By the continuity of $g_\Omega(\cdot,\cdot)$ on $\Omega\times \Omega\setminus \diag$, 
the functions
$\phi h_k$ converge pointwise to $\phi g_\Omega(x,\cdot)$ on $\Omega$,
and using dominated convergence this convergence also holds in $L^2(\Omega)$.

For the remainder of the proof we identify $L^2(\Omega)$ with the subspace \[ \Sett{f \in L^2(\sX)}{f= 0,\hspace{3mm} \mu-\mbox{a.e. on } \Omega^\compl}.\] Similarly, we view $\sF^0(\Omega)$ as a subspace of $\sF(\Omega)$ --see \cite[ (3.2.2) and Theorem 3.3.9]{CF}. In particular, we can view the 
functions $\phi h_i$ as functions over $\sX$. By \cite[Theorem 1.4.2(ii),(iii)]{FOT},
$\phi h_i = \phi (h_i \wedge M) \in \sF^0(\Omega)$, and $\phi^2 h_i = \phi^2 (h_i \wedge M) \in \sF^0(\Omega)$. 

We now show that $\phi h_i$ is Cauchy in the seminorm induced by $\sE(\cdot,\cdot)$. 
By the Leibniz rule (cf. \cite[Lemma 3.2.5]{FOT}) we have
\begin{equation}\label{e:gn7}
\sE(\phi(h_i-h_j),\phi (h_i - h_j)) = \int_{\sX} (h_i - h_j)^2 \, d\Gamma(\phi,\phi) + \sE(h_i-h_j, \phi^2(h_i-h_j)).
\end{equation}
Since $\one_{B_i} - \one_{B_j}$ and $\vp^2(h_i-h_j)$ have disjoint support
the second term is zero by Lemma \ref{l:5.1}(ii).

For the first term in \eqref{e:gn7}, we use the fact that the function
$h_i$ vanishes on 
$\Omega^\compl$ together with strong locality to obtain
\begin{equation}\label{e:gn8}
\int_{\sX} (h_i - h_j)^2 \, d\Gamma(\phi,\phi) = \int_{\overline{B(x,3r)} \setminus B(x,2r)} (h_i - h_j)^2 \, d\Gamma(\phi,\phi).
\end{equation}
Let $F = \overline{B(x,3r)} \setminus B(x,2r)$. Since
$\overline{B(x,r)} \times F$ is compact, by Assumption \ref{a-green}(iii)
the function $g_\Omega(\cdot,\cdot)$ is uniformly continuous on
$\overline{B(x,r)} \times F$.
This in turn implies that $h_i$ converges uniformly to $g_\Omega(x,\cdot)$ 
as $i \to \infty$ on $F$, and so by \eqref{e:gn8} 
we have that $\phi h_i$ is Cauchy in the $\left(\sE_\Omega(\cdot,\cdot) \right)^{1/2}$-seminorm.
Since $\phi h_i$ converges pointwise and in $L^2(\Omega)$ to $\phi g_\Omega(x,\cdot)$,
and $(\sE_\Omega,\sF^0(\Omega))$ is a closed form, this implies that $\phi g_\Omega(x,\cdot) \in \sF^0(\Omega)$.

Finally, we show that $g_\Omega(x,\cdot)$ is harmonic on $\Omega\setminus \set{x}$. 
Let $\psi \in \sF_c(\Omega \setminus \set{x})$, and let 
$V \Subset \Omega$ be a precompact open set containing $\supp(\psi)$ such that $d(x,V)>0$. 
Choose $r>0$ such that $B(x,4r) \cap V = \emptyset$, and let $\vp$ and $h_k$ be 
as defined above. Then as $\vp \equiv 1$ on $V$, using strong locality,
\[ 
\sE(\vp h_k, \psi) = \sE( h_k, \psi) = \mu(B_k)^{-1} \langle \one_{B_k},\psi \rangle =0.
\]
As $\vp h_k$ converge to $\phi g_\Omega(x,\cdot)$ in the $\sE_1(\cdot, \cdot)^{1/2}$ norm,
it follows that $\sE(\phi g_\Omega(x,\cdot), \psi )=0$. 
%
This allows us to conclude that $g_\Omega(x,\cdot)$ is harmonic on  $\Omega\setminus \set{x}$.  
\qed
 
The elliptic Harnack inequality 
enables us to relate capacity and Green functions, and also
to control  their fluctuations  on bounded regions of $\sX$.

\begin{proposition}\label{p-consq} (See \cite[Section 3]{BM}).
Let $(\sX,d,\mu,\sE,\sF)$  be a metric measure Dirichlet space satisfying the EHI and 
Assumption \ref{a-green}. 
Then the following hold:
\begin{enumerate}[(a)]
	 \item 	 For all $A_1,A_2 \in (1, \infty)$, there exists $C_0= C_0(A_1,A_2,C_H)>1$ such that for all bounded open sets $D$  
and for all $x_0 \in \sX, r>0$ that satisfy $B(x_0,A_1r) \subset D$, we have
	 \[
	 g_D(x_1,y_1) \le C_0 g_D(x_2,y_2) \q \forall x_1,y_1,x_2,y_2 \in  B(x_0,r),
	\]
	satisfying $d(x_i,y_i) \ge r/A_2$, for $i=1,2$.
	\item For all $A \in (1, \infty)$, there exists $C_1= C_1(A,C_H)>1$ such that for all bounded open sets $D$  
and for all $x_0 \in \sX, r>0$ that satisfy $B(x_0,A r) \subset D$, we have
	\[
	\inf_{y \in \partial B(x_0,r)} g_D(x_0,y)\le \Cpc_D\left( \overline{B(x_0,r)}\right)^{-1} \le  \Cpc_D\left({B(x_0,r)}\right)^{-1}  \le C_1 \inf_{y \in \partial B(x_0,r)} g_D(x_0,y).
	\]
	\item For all $1\le A_1 \le A_2 < \infty$ and $a \in (0,1]$ there exists
	$C_2= C_2(a,A_1,A_2,C_H)>1$ such that for $x \in \sX$, and $r>0$  with
	 $r \le \operatorname{diameter}(\sX)/5 A$,
	\[
	   \Cpc_{B(x_0,A_2r)}\left(B(x_0,ar)\right) 
	    \le  \Cpc_{B(x_0,A_1r)}\left(B(x_0,r)\right) \le  C_2 \Cpc_{B(x_0,A_2r)}\left(B(x_0,a r)\right). 
 	\]
	\item
	 For all $A_2 > A_1 \ge 2$ there exists $C_3= C_3(A_1,A_2,C_H)>1$ such that for 
	 all $x,y \in \sX$, with $d(x,y)=r>0$ and 
	 such that $r \le \operatorname{diameter}(\sX)/5 A_2$,
	\[
	g_{B(x,A_1 r)}( x,y) \le g_{B(x,A_2r)}( x,y) \le C_3 g_{B(x,A_1 r)}( x,y). 
 	\]
\item $(\sX,d)$ satisfies metric doubling. 
	\end{enumerate}
\end{proposition}

The statements given above are slightly stronger than those in
\cite[Section 3]{BM}, 
but Proposition \ref{p-consq} easily follows  from the results
there using additional chaining arguments.
%

We will need the following maximum principle for Green functions. 

\begin{lemma}\label{L:maxp}
Suppose that $(\sX,d,\mu,\sE,\sF)$ and $\Omega$ satisfy Assumption \ref{a-green}.
Let $c_0 \in (0,1)$, $y,y^* \in \Omega$, such that $B(y^*,r) \Subset \Omega$ and
 \[
 g_\Omega(y,x) \ge c_0 g_\Omega(y^*,x) \q \text{ for all } x \in \partial B(y^*,r).
 \]
 Then 
  \[
 g_\Omega(y,x) \ge c_0 g_\Omega(y^*,x) \q \text{ for all } 
 x \in \Omega \setminus ( \{y\} \cup B(y^*,r)) .
 \]
 \end{lemma}

 \proof If $y \in B(y^*,r)$, then the function
 $g_\Omega(y, \cdot) - c_0 g_\Omega(y^*, \cdot) $ is bounded and harmonic
 on $\Omega \setminus  B(y^*,r)$, so the result follows by the maximum principle
 in \cite[Lemma 4.1(ii)]{GH}.
 
 Now suppose that $y \not\in B(y^*,r)$. 
 Choose $r'>0$ small enough so that $B(y,4 r') \subset  \Omega \setminus B(y^*,r)$,
 and as in \eqref{e:appxG} set
 \[
  h_n(x) = \mu(B(y, r'/n))^{-1} \int_{B(y, r'/n)} g_\Omega(z,x) \mu(dz). 
 \]
 Then $h_n \in \sF^0(\Omega)$ and 
 $h_n \to g_\Omega(y, \cdot)$ pointwise on $\Omega\setminus \{y\}$.
Let $M= 2 \sup_{ z \in \partial B(y^*,r)} g_\Omega(z,y^*)$.
By \cite[Corollary 2.2.2 and Lemma 2.2.10]{FOT} 
the functions $M \wedge h_n$ are bounded and superharmonic on $\Omega$.

Set $f_n = M \wedge h_n - c_1 g_\Omega(y^*, \cdot)$, where  $c_1 \in (0,c_0)$.  
Then $f_n$ is a bounded, 
superharmonic function in $\Omega \setminus B(y^*,r)$ that is non-negative on the 
 boundary of $\Omega \setminus B(y^*,r)$ for all sufficiently large $n$. By the maximum principle in 
 \cite[Lemma 4.1(ii)]{GH} we obtain that $f_n \ge 0$ in 
 $\Omega \setminus B(y^*,r)$ for all sufficiently large $n$.
 Since $c_1 \in (0,c_0)$ was arbitrary, we obtain the desired conclusion by letting $n \to \infty$.
 \qed


\begin{remark} \label{R:com}
{\rm Suppose that Assumption \ref{a-green} holds for $(\sX,d,\mu,\sE,\sF)$.
Let $\mu'$ be a measure which is mutually absolutely continuous
with respect to the measure $\mu$, and suppose
that $d\mu'/d\mu$ is uniformly bounded away from 0 and infinity
on bounded sets. It
is straightforward to verify that if $\Omega$ is a bounded domain
and the operator $G'_\Omega$ is defined by
$$ G'_\Omega f(x) = \int_{\Omega} g_\Omega(x,y)f(y) \mu'(dy), 
\mbox{ for } f \in L^2(\Omega, \mu'), $$
then $G'_\Omega$ is the Green operator on the domain $\Omega$ for the
Dirichlet form $(\sE,\sF')$ on $L^2(\sX,\mu')$, where $\sF'$ is the domain of the time-changed Dirichlet space (Cf. \cite[p. 275]{FOT}).
It follows that $g_\Omega(x,y)$ is the Green function for both 
$(\sX,d,\mu,\sE,\sF)$ and $(\sX,d,\mu',\sE,\sF')$, and thus that
Assumption \ref{a-green} holds for $(\sX,d,\mu',\sE,\sF')$.
}\end{remark}

\subsection{Sufficient conditions for Assumption \ref{a-green}}

We begin by recalling the definition of an ultracontractive semigroup, 
a notion introduced in \cite{DS}.

\begin{definition}\label{d-ultra} 
{\rm Let $(\sX,d,\mu)$ be a metric measure space.
Let $(P_t)_{t \ge 0}$ be the semigroup associated with the Dirichlet form $(\sE,\sF)$ 
on $L^2(\sX,\mu)$. We say that the semigroup $(P_t)_{t \ge 0}$ is \emph{ultracontractive} if 
$P_t$ is a bounded operator from $L^2(\sX,\mu)$ to $L^\infty(\sX,\mu)$ for all $t > 0$. 
We say that $(\sE,\sF)$ 
is \emph{ultracontractive} if the 
associated heat semigroup is ultracontractive.
}\end{definition}

We will use the weaker notion of local ultracontractivity introduced 
in \cite[Definition 2.11]{GT12}.

\begin{definition} \label{d-localultra}
{\rm	We say that a MMD space $(\sX,d,\mu,\sE,F)$ is {\em locally ultracontractive} 
if for all open balls $B$, the killed heat semigroup $(P_t^B)$ 
given by Definition \ref{d-dirdir} is ultracontractive.
}\end{definition}

It is well-known that ultracontractivity of a semigroup is equivalent to the existence of an essentially bounded heat kernel at all strictly positive times. 

\begin{lemma} (\cite[Lemma 2.1.2]{Dav})
Let $U$ be a bounded open subset of $\sX$. \\
(a) Suppose that $(P^U_t)$ is ultracontractive.
Then for each $t>0$ the operator
$P_t^U$ has an integral kernel $p^U(t,\cdot,\cdot)$ which is
jointly measurable in $U \times U$ and satisfies
\[
0 \le p^U(t,x,y) \le \norm{P_{t/2}^U}_{L^2(\mu) \to L^\infty(\mu)}^2
\mbox{ for $\mu \times \mu$-a.e.~$(x,y) \in U \times U$. }
\]
(b) If $P_t^U$ has an integral kernel $p^U(t,x,y)$  satisfying
\[
0 \le  p^U(t,x,y) \le a_t < \infty
\]
for all $t >0$ and for $\mu \times \mu$-a.e.~$(x,y) \in U \times U$, 
then $(P_t^U)_{t \ge 0}$  is ultracontractive with
\[
\norm{P_{t}^U}_{L^2(\mu) \to L^\infty(\mu)} \le a_t^{1/2} \text{ for all $t > 0$. }
\]
\end{lemma}

The issue of joint measurability is clarified in \cite[p. 1227]{GT12}.

We now introduce a second assumption on the Dirichlet form $(\sE,\sF)$,
and will prove below that it implies Assumption \ref{a-green}.

\begin{assumption}\label{A:locult}
{\rm Let $(\sX,d)$ be a complete, locally compact, separable, length space and 
let $\mu$ be a non-atomic Radon measure on $(\sX,d)$ with full support.
Let $(\sE,\sF)$ be a strongly local, regular, Dirichlet form on $L^2(\sX,\mu)$. 
We assume that  the MMD space $(\sX,d,\mu,\sE,\sF)$ is locally ultracontractive 
and that $\lambda_{\min}(\Omega)  > 0$
for any non-empty bounded open set $\Omega \subset \sX$ with 
$\operatorname{diameter}(\Omega,d) \le \operatorname{diameter}(\sX,d)/5$.
}\end{assumption}




This assumption gives
the existence of a Green function satisfying Assumption  \ref{a-green}.

\begin{lemma}\label{l:green}( See \cite[Lemma 5.2 and 5.3]{GH})).
Let $(\sX,d,\mu,\sE,\sF)$  be a metric measure Dirichlet space which satisfies
the EHI$_{\rm loc}$ and Assumption \ref{A:locult}. 
Then Assumption \ref{a-green} holds.
\end{lemma}

\begin{remark} \label{R:gh-error}
{\rm A similar result  is stated in \cite[Lemma 5.2]{GH}. 
Unfortunately, the proof in \cite{GH} of Assumption \ref{a-green}(i) and (iv) have 
gaps, which we
do not know how to fix without the extra hypothesis of local ultracontractivity. 
For (i) the problem occurs in the proof of \cite[(5.8)]{GH} from \cite[(5.7)]{GH}. 
In particular, while one has in the notation of \cite{GH} that
$G^\Omega f_k \to G^\Omega f$ in $L^2(\Omega)$, this does not imply
pointwise convergence. On the other hand the proof of \cite[(5.8)]{GH} 
does require pointwise convergence at the specific point $x$.
	
The following example helps to illustrate this gap. Consider
the Dirichlet form $\sE(f,f)=\norm{f}_2^2$ on $\bR^n$; this
satisfies \cite[(5.7)]{GH} for any bounded open domain $\Omega$ with $g_x^\Omega \equiv 0$ 
but it fails to satisfy \cite[(5.8)]{GH}. This Dirichlet form 
does not satisfy the hypothesis of \cite[Lemma 5.2]{GH} since it is local 
rather than strongly local, but it still  illustrates the problem,  since strong locality 
was not used in the proof of \cite[(5.8)]{GH} from \cite[(5.7)]{GH}.
} \end{remark}

\noindent {\em Proof of Lemma \ref{l:green}.}
Let $\Omega$ be a non-empty bounded open set with 
$\operatorname{diameter}(\Omega,d) \le \operatorname{diameter}(\sX,d)/5$;
we need to verify properties (i)--(iv) of Assumption \ref{a-green}.

We use the construction in \cite{GH}. We denote by $g_\Omega(\cdot,\cdot)$ the function 
constructed in \cite[Lemma 5.2]{GH} off the diagonal, 
and extend it to $\Omega\times \Omega$ by taking $g_\Omega$ 
equal to $0$ on the diagonal. 
By \cite[Lemma 5.2]{GH} the function $g_\Omega(\cdot,\cdot)$ satisfies (ii) and (iii).
(The proofs of (ii) and (iii) do not use \cite[(5.8)]{GH}.)

Next, we show (i), using the additional hypothesis of local ultracontractivity.
Define the operators 
\[S_t= P^\Omega_t \circ G^\Omega=   G^\Omega \circ P^\Omega_t,\q t  \ge 0. \]
Formally we have
$S_t = \int_t^\infty P_s ds.$
Since $P^\Omega_t$ is a contraction on all $L^p(\Omega)$, we have by Lemma \ref{l:5.1}(i)
\begin{equation}
\norm{S_t}_{L^2 \mapsto L^2} \le \norm{G^\Omega}_{L^2 \mapsto L^2} \norm{P_t^{\Omega}}_{L^2 \mapsto L^2} < \infty. \label{e:gn-1}
\end{equation}
 Therefore by \cite[Lemma 1.4.1]{FOT}, for each 
 $t \ge 0$ there exists a positive symmetric Radon measure $\sigma_t$ on 
 $\Omega \times \Omega$ such that for all functions $f_1,f_2 \in L^2(\Omega)$, we have
 \begin{equation}\label{e:gn0}
 \langle f_1, S_t f_2 \rangle = \int_{\Omega \times \Omega} f_1(x) f_2(y) \sigma_t( dx,dy).
 \end{equation}
By \cite[Lemma 1.4.1]{FOT}
$S_{t+r}-S_t$ is a positive symmetric operator on $L^2(\Omega)$. Thus
$(\sigma_t, t \in \bR_+)$ is a  family of symmetric positive measures 
on $\Omega \times \Omega$ with $\sigma_s \ge \sigma_t$ for all $0\le s \le t$.
Note that the measures $\sigma_t$ are finite, since for any $t\ge 0$
by using Lemma \ref{l:5.1}(i), we have
\be  \label{e:gn1}
 \sigma_t(\Omega \times \Omega)  \le \sigma_0(\Omega \times \Omega)
 = \langle \one_\Omega, G^\Omega \one_\Omega \rangle
 \le ||  \one_\Omega||_2^2 \lambda_{\min}(\Omega)^{-1} = \mu(\Omega) \lambda_{\min}(\Omega)^{-1}.
\ee

Let $A,B$ be measurable subsets of $\Omega$.
Since $P_t^\Omega$ is a strongly continuous semigroup, we have
\begin{equation*}
\sigma_0(A \times B) =  \langle G^\Omega \one_A , \one_B \rangle= \lim_{t \downarrow 0}  \langle S_t \one_A , \one_B \rangle = \lim_{t \downarrow 0} \sigma_t(A \times B).
\end{equation*}
The above equation implies that, for all measurable subsets 
$F \subset \Omega \times \Omega$, we have
\begin{equation}\label{e:gn2}
\lim_{n \to \infty} \sigma_{1/n}( F) = \sigma_0 (F).
\end{equation}

For each $t >0$ by local ultracontractivity, we have
\begin{equation} \label{e:gn3}
\norm{S_t}_{L^2 \to L^\infty} 
\le  \norm{G^\Omega}_{L^2 \to L^2}\norm{P^\Omega_{t}}_{L^2 \to L^\infty} < \infty.
\end{equation}
By \cite[Lemma 3.3]{GH2}, there exists a jointly measurable function $s_t(\cdot,\cdot)$ on $\Omega \times \Omega$ such that
 \begin{equation}\label{e:gn4}
 \langle f_1, S_t f_2 \rangle_{L^2(\Omega)} = \int_{\Omega} \int_{\Omega} f_1(x) f_2(y) s_t(x,y) \mu(dx) \mu(dy),
 \end{equation}
 for all $f_1,f_2 \in L^2(\Omega)$.  By \eqref{e:gn0} and \eqref{e:gn4}, we have
 \be \label{e:gn4a}
 \sigma_{1/n}(dx,dy) = s_{1/n}(x,y) \mu(dx) \mu(dy)
 \ee
for all $n \in \bN$. Therefore by \eqref{e:gn4a},
\eqref{e:gn1}, \eqref{e:gn2} and Vitali-Hahn-Saks theorem (cf. \cite[p. 70]{Yos}), the measure $\sigma_0$ is absolutely continuous with respect to the product measure $\mu \times \mu$ on $\Omega \times \Omega$. Let  $s(\cdot,\cdot)$ be the Radon-Nikoym derivative of $\sigma_0$ with respect to $\mu \times \mu$. By \eqref{e:gn0} and Fubini's theorem, for all $f \in L^2(\Omega)$ and for almost all $x \in \Omega$,
	\begin{equation}\label{e:gn4p}
	G^\Omega f(x) = \int_\Omega s(x,y) f(y) \, \mu(dy).
	\end{equation}
If $B$ and $B'$ are disjoint open balls in $\Omega$, then for all $f_1 \in L^2(B), f_2 \in L^2(B')$,  we have
\begin{align} \nn
\langle G^\Omega f_1, f_2 \rangle &= \int_{B} \int_{B'} f_1(x) f_2(y) s(x,y)\mu(dy) \mu(dx) \\ 
\label{e:gn4z}
&=  \int_{B} \int_{B'} f_1(x) f_2(y) g_\Omega(x,y)\mu(dy) \mu(dx).
\end{align}
We used  \cite[(5.7)]{GH}, along with \eqref{e:gn4}, to obtain the above equation.
By the same argument as in \cite[Lemma 3.6(a)]{GH2}, \eqref{e:gn4z} implies that
\begin{equation}\label{e:gn5}
s(x,y)=g_\Omega(x,y) \mbox { for $\mu \times \mu$-almost every $(x,y) \in B\times B'$. }
\end{equation}

By an easy covering argument 
$\Omega \times \Omega \setminus \diag$ can be covered by countably 
many sets of the form $B_i \times B_i', i \in \bN$, such that $B_i$ and $B_i'$ are 
disjoint balls contained in $\Omega$. Therefore by \eqref{e:gn5}, we have
\begin{equation}\label{e:gn6}
s(x,y)=g_\Omega(x,y)
\end{equation}
for almost every $(x,y) \in \Omega \times \Omega$. In the last line we used that 
since $\mu$ is non-atomic the  diagonal has measure zero.  
Combining \eqref{e:gn6} and \eqref{e:gn4p} gives property (i). 

The maximum principles in (iv) are proved in \cite[Lemma 5.3]{GH} by showing the corresponding maximum principles for the approximations of Green functions 
given by \eqref{e:appxG}. This maximum principle for the
Green functions follows from 
\cite[Lemma 4.1]{GH}, 
provided  that the functions $h_k$ in \eqref{e:appxG} satisfy the three
properties that $h_k \in \sF^0(\Omega)$, $h_k$ is superharmonic, 
and $h_k \in L^\infty(\Omega)$. As in \cite{GH}, the first two conditions can be 
checked by using Lemma \ref{l:5.1}(i) and (ii). 
To verify that $h_k \in L^\infty(\Omega)$ it is
sufficient to prove that $G^\Omega \one_\Omega \in L^\infty(\Omega)$. By Lemma \ref{l:5.1}(iii) and the semigroup property, for any $t>0$,
\[
G^\Omega \one_\Omega = \int_0^t P_s^\Omega \one_\Omega \,ds + P_t^\Omega \circ G^\Omega \one_\Omega.
\]
For the first term, we use $\norm{P_s^\Omega}_{L^\infty \to L^\infty} \le 1$, and for the second term we use local ultracontractivity 
and Lemma \ref{l:5.1}(i).
\qed

We now introduce some conditions which imply local ultracontractivity.
We say that a function $u=u(x,t)$ is {\em caloric} in a region 
$Q \subset \sX \times (0, \infty)$ if $u$ is a weak solution of $(\partial_t+ \sL)u=0$ in $Q$;
here $\sL$ is the generator corresponding to the  Dirichlet form $(\sE,\sF, L^2(\sX,\mu))$.
Let $\Psi:[0,\infty) \to [0,\infty)$ have the property that there exist
constants $1< \beta_1 \le \beta_2<\infty$ and $C>0$ such that
\be \label{e:Psi_cond}
 C^{-1} (R/r)^{\beta_1} \le \Psi(R)/\Psi(r) \le  C (R/r)^{\beta_2}. 
\ee

\begin{definition} \label{d-locvdpiphi}
We say that a MMD space $(\sX,d,\mu,\sE,\sF)$ satisfies the local volume doubling 
property $(\operatorname{VD})_{\operatorname{loc}}$, if there exist $R \in (0,\infty]$,
$C_{VD}>0$ such that
\begin{equation*} \label{vd-loc}
V(x,2 r)  \le C_{VD} V(x,r)\tag*{$(\operatorname{VD})_{\operatorname{loc}}$}
\mbox{ for all $x \in \sX$ and for all $0<r \le R$.}
\end{equation*}

We say a MMD space $(\sX,d,\mu,\sE,\sF)$ satisfies the local Poincar\'e inequality 
$(\operatorname{PI(\Psi)})_{\operatorname{loc}}$, if there exist $R \in (0,\infty]$,
$C_{PI}>0$, and $A\ge 1$ such that
\begin{equation*} \label{pi2-loc}
  \int_{B(x,r)} \abs{f-f_{B(x,r)}}^2 \,d\mu \le C_{PI} \Psi(r) \int_{B(x,Ar)} d\Gamma(f,f) \tag*{$(\operatorname{PI(\Psi)})_{\operatorname{loc}}$}
	\end{equation*}
for all $x \in \sX$ and for all $0<r \le R$, where $\Gamma(f,f)$ denotes the energy measure, and $f_{B(x,r)}= \frac{1}{\mu(B(x,r))}\int_{B(x,r)} f \, d\mu$. \\

We say that a MMD space  $(\sX,d,\mu,\sE,\sF)$ satisfies the local parabolic Harnack 
inequality $(\operatorname{PHI(\Psi)})_{\operatorname{loc}}$, if there exist 
$R \in (0,\infty]$, $C_{PHI}>0$
such that for all $x \in \sX$, for all $0< r \le R$, any non-negative 
caloric function $u$ on $(0,r^2) \times B(x,r)$ satisfies
	\begin{equation*} \label{phi2-loc}
	\sup_{(\Psi(r)/4,\Psi(r)/2) \times B(x,r/2)} u \le C_{PHI} \inf_{(3\Psi(r)/4,\Psi(r)) \times B(x,r/2)} u, \tag*{$(\operatorname{PHI(\Psi)})_{\operatorname{loc}}$}
	\end{equation*} 
We will write 
$(\operatorname{PI(\beta)})_{\operatorname{loc}}$ and
$(\operatorname{PHI(\beta)})_{\operatorname{loc}}$ for the
conditions $(\operatorname{PI(\Psi)})_{\operatorname{loc}}$
and $(\operatorname{PHI(\Psi)})_{\operatorname{loc}}$ if $\Psi(r)=r^\beta$.
\end{definition}

\begin{lemma}\label{l-ex-adir}
Let $(\sX,d)$ be a complete, locally compact, separable, length space with 
$\operatorname{diameter}(\sX)=\infty$,  let $\mu$ be a Radon measure on $(\sX,d)$ 
with full support and  let $(\sE,\sF)$ be a strongly local, regular, Dirichlet form on 
$L^2(\sX,\mu)$. 
If $(\sX,d,\mu,\sE,\sF)$  satisfies the properties
 $(\operatorname{VD})_{\operatorname{loc}}$ and 
 $(\operatorname{PI(2)})_{\operatorname{loc}}$,
then $(\sX,d,\mu,\sE,\sF)$  satisfies Assumption \ref{a-green}, and the
property {\rm EHI}$_{\rm loc}$.
 \end{lemma}

\proof
First, the property $(\operatorname{PHI(2)})_{\operatorname{loc}}$
is satisfied;
this is immediate from \cite[Theorem 8.1]{CS} and \cite[Theorem 2.7]{HS},
which prove that the property $(\operatorname{PHI(2)})_{\operatorname{loc}}$
is equivalent to the conjunction of the properties
$(\operatorname{PI(2)})_{\operatorname{loc}}$ and 
$(\operatorname{VD})_{\operatorname{loc}}$. 

To prove Assumption \ref{a-green}
it is sufficient to verify  
the property {\rm EHI}$_{\rm loc}$ and Assumption \ref{A:locult}. 
Of these, the property {\rm EHI}$_{\rm loc}$ follows immediately from the local PHI.

The heat kernel  corresponding to 
$(\sX,d,\mu,\sE,\sF)$  satisfies Gaussian upper bounds for small times by
\cite[Theorem 2.7]{HS}. Since the  heat kernel of the killed semigroup
is dominated by the heat kernel of $(\sX,d,\mu,\sE,\sF)$, 
local ultracontractivity follows using the property \ref{vd-loc}.
The fact that $\mu$ is non-atomic follows from the property \ref{vd-loc} 
due to a reverse volume doubling property -- see \cite[(2.5)]{HS}. 

By domain monotonicity, it suffices to verify that $\lambda_{\min}(B(x,r))>0$ for all 
balls $B=B(x,r)$ with $0<r < \operatorname{diameter}(\sX)/4$. 
Consider a ball $B(z,r)$ such that $B(x,r)\cap B(z,r)=\emptyset$ and $d(x,z)\le 3r$. 
By the Gaussian lower bound for small times \cite[Theorem 2.7]{HS} and 
the property \ref{vd-loc}, there exists $t_0>0,\delta \in (0,1)$ such that
\[
\bP^y (X_{t_0}  \in B(z,r)) \ge \delta, \q \forall y \in B(x,r),
\]
where $(X_t)_{t \ge 0}$ is the diffusion corresponding to the MMD space $(\sX,d,\mu,\sE,\sF)$. 
This implies that $P_{t_0}^B \one_B \le (1- \delta) \one_B$, which in turn implies that
\[
P_{t}^B \one_B \le (1- \delta)^{ \lfloor(t/t_0) \rfloor} \one_B, \q \forall t \ge 0.
\]
It follows that $\norm{G^B \one_B}_{L^\infty} < \infty$. 
By Riesz--Thorin interpolation we have 
$\norm{G^B}_{L^2 \to L^2}  \le \norm{G^B}_{L^1 \to L^1}^{1/2} \norm{G^B}_{L^\infty \to L^\infty} ^{1/2}$, while by duality 
$\norm{G^B}_{L^\infty \to L^\infty} = \norm{G^B}_{L^1 \to L^1}$. Thus
$$
\lambda_{\min}(B)^{-1}= \norm{G^B}_{L^2 \to L^2} \le \norm{G^B}_{L^\infty \to L^\infty} = \norm{G^B \one_B}_{L^\infty} < \infty. $$
\qed

\subsection{Examples}

In this section, we give some examples of MMD spaces which satisfy 
Assumption \ref{a-green}: 
weighted Riemannian manifolds and cable systems of weighted graphs.
We also briefly describe some  
classes of regular fractals which satisfy Assumption \ref{a-green}
-- see Remark \ref{R:Lier}.

\sm{\bf Example 1.} (Weighted Riemannian manifolds.)
Let $(\sM,g)$ be a Riemannian manifold,
and $\nu$ and  $\nabla$  denote the Riemannian measure and the 
Riemannian gradient respectively. Write $d=d_g$ for the Riemannian distance function.
A {\em weighted manifold} $(\sM,g,\mu)$ is a Riemmanian manifold $(\sM,g)$ 
endowed with a measure $\mu$ that has a smooth (strictly) positive density $w$ with 
respect to the Riemannian measure $\nu$. 
The {\em weighted Laplace operator} $\Delta_\mu$ on $(\sM,g,\mu)$
is given  by
\[
\Delta_\mu f= \Delta f + g\left( \nabla \left( \ln w \right), \nabla f  \right), \quad f \in \sC^\infty(\sM).
\]
We say that the weighted manifold $(\sM,g, \mu)$ has {\em controlled weights} 
if $w$ satisfies
\[
\sup_{x,y \in \sM : d(x,y) \le 1} \frac{w(x)}{w(y)} < \infty.
\]
The construction of heat kernel, Markov semigroup and Brownian motion for a 
weighted Riemannian manifold $(\sM,g,\mu)$ is outlined in \cite[Sections 3 and 8]{Gri06}. 
The corresponding Dirichlet form on $L^2(\sM,\mu)$ given by
\[
\sE_w(f_1,f_2)= \int_{\sX} g(\nabla f_1, \nabla f_2)\,d\mu, \q f_1,f_2 \in \sF,
\] 
where $\sF$ is the weighted Sobolev space of functions in $L^2(\sM,\mu)$ whose distributional gradient is also in $L^2(\sM,\mu)$.  See \cite{Gri06} and \cite[pp. 75--76]{CF}
for more details. 

\sm {\bf Example 2.} (Cable systems of weighted graphs.) 
Let $\bG = (\bV,E)$ be an infinite graph, such that each vertex 
$x$ has finite degree. For $x \in \bV$ we
write $x \sim y$ if $\{x,y\} \in E$. 
Let $w: E \to (0,\infty)$ be a function
which assigns weight $w_e$ to the edge $e$. We write
$w_{xy}$ for $w_{\{x,y\}}$, and define
$$ w_x = \sum_{y \sim x} w_{xy}. $$
We call $(\bV,E,w)$ a {\em weighted graph}. An {\em unweighted graph} has
$w_e = 1$ for all $e\in E$. We say that $\bG$ has {\em controlled weights} if there exists
$p_0>0$ such that
\be
  \frac{w_{xy}}{w_x} \ge p_0 \, \hbox{ for all } x\in \bV, y \sim x.
\ee

The {\em cable system of a weighted graph}
gives a natural embedding of a graph in a connected metric length space.
Choose a direction for each edge $e \in E$, let $(I_e, e \in E)$ be a collection
of copies of the open unit interval, and set
$$ \sX = \bV \cup \bigcup_{e \in E} I_e. $$
(We call the sets $I_e$ {\em cables}).
We define a metric $d_c$ on $\sX$ by using Euclidean distance on each cable,
and then extending it to a metric
on $\sX$; note that this agrees with the graph metric for $x,y \in \bV$.
Let $m$ be the measure on $\sX$ which assigns zero mass to 
points in $\bV$, and mass $w_e |s-t|$ to any interval $(s,t) \subset I_e$.
It is straightforward to check that $(\sX,d_c,m)$ is a MMD space.
For more details on this construction see \cite{V,BB3}.

We say that a function $f$ on $\sX$ is piecewise differentiable
if it is continuous at each vertex $x \in \bV$, is differentiable on each cable, and
has one sided derivatives at the endpoints. Let $\sF_{d}$ be the set
of piecewise differentiable functions $f$ with compact support.
Given two such functions we set
\begin{align*}
d\Gam(f,g)(t) = f'(t)g'(t) m(dt), \q
\sE(f,g) = \int_{\sX}  d\Gam(f,g)(t), \q f,g \in \sF_d,
\end{align*}
and let $\sF$ be the completion of $\sF_d$ with respect to the 
$\sE_1^{1/2}$norm.
We extend $\sE$ to $\sF$; it is straightforward to verify 
that $(\sE,\sF)$ is a closed regular strongly local Dirichlet form. 
We call $(\sX,d_c,m,\sE,\sF)$ the {\em cable system} of the graph $\bG$.

%

We will now show that both these examples satisfy 
the conditions  $(\operatorname{VD})_{\operatorname{loc}}$ and 
$(\operatorname{PI(2)})_{\operatorname{loc}}$, 
and therefore Assumption \ref{a-green}.

%
%


\begin{lemma} \label{l-ex-philoc}
(a) Let $(\sM,g,\mu)$ be a weighted Riemmanian manifold with controlled weights $w$ 
which is quasi isometric to a Riemannian manifold $(\sM',g')$ with Ricci curvature bounded 
below. Then the MMD space $(\sM,d_g,\mu,\sE_w)$ satisfies the conditions
 $(\operatorname{VD})_{\operatorname{loc}}$ and 
 $(\operatorname{PI(2)})_{\operatorname{loc}}$. \\
(b) Let $\bG$ be a weighted graph with controlled weights. Then the corresponding cable system satisfies the conditions
 $(\operatorname{VD})_{\operatorname{loc}}$ and 
 $(\operatorname{PI(2)})_{\operatorname{loc}}$. 
\end{lemma}

\proof 
(a) The properties \ref{vd-loc} and 
$(\operatorname{PI(2)})_{\operatorname{loc}}$
for $(\sM',g')$ 
follow from the Bishop-Gromov volume comparison theorem \cite[Theorem III.4.5]{Cha} and Buser's Poincar\'{e} inequality (see \cite[Lemma 5.3.2]{Sal02}) respectively. 
Since quasi isometry only changes distances and volumes by at most a constant factor, 
we have that \ref{vd-loc} and $(\operatorname{PI(2)})_{\operatorname{loc}}$
also hold for $(\sM,g)$. The controlled weights
condition on $w$ implies that these two conditions also hold for $(\sM,g,\mu)$. \\
(b) Using the controlled weights condition and the uniform bound on vertex degree, 
one can easily obtain  the two properties \ref{vd-loc} and 
$(\operatorname{PI(2)})_{\operatorname{loc}}$. \qed

\begin{remark} \label{R:Lier}
{\rm  
The paper \cite{Lie} proves a BHP on MMD spaces which 
are length spaces and 
satisfy a weak heat kernel estimate associated with 
a space time scaling function $\Psi$, where
$\Psi$ satisfies the condition \eqref{e:Psi_cond}.
By \cite[Theorem 3.2]{BGK} these spaces 
satisfy $(\operatorname{PHI(\Psi)})_{\operatorname{loc}}$ with $R=\infty$. 
The same argument as in Lemma \ref{l-ex-adir} then proves that
these spaces satisfy Assumption \ref{a-green}.

Examples of spaces of this type are the Sierpinski gasket, nested
fractals, and generalized Sierpinski carpets -- see 
\cite{BP,Kum1,BB}.
%
}\end{remark}
%


\section{Proof of Boundary Harnack Principle} \label{sec:bhp}

In this section we prove Theorem \ref{T:BHP}. 
For the remainder of the section, we assume the hypotheses of Theorem \ref{T:BHP},
and will fix a $(c_U, C_U)$-inner uniform domain $U$. 
We can assume that $c_U\le \half \le 2 \le C_U$, and will
also assume that the EHI holds with constants $\delta=\half$ and $C_H$.
We will use $A_i$ to denote constants which just depend 
on the constants $c_U$ and $C_U$; other constants will
depend  on $c_U, C_U$ and $C_H$.

Since by  Proposition \ref{p-consq}(e), $(\sX,d)$ has the metric doubling property,
we can use Lemma \ref{l:ballcomp}. In addition we will assume that
$$ \operatorname{diameter}(U)=\infty, $$
so that $R(U)=\infty$. The proof of the general case is the same except
that we need to ensure that the balls $B_U(\xi,s)$ considered in the argument
are all small enough so that they do not equal $U$.

\begin{definition}[Capacitary width]
{\rm For an open set $V \subset \sX$ and $\eta \in (0,1)$, define 
the capacitary width $w_\eta(V)$ by
\be \label{e:w-eta}
w_\eta(V)= \inf\Sett{r>0}{ \frac{\Cpc_{B(x,2r)} \left( \overline{B(x,r)} \setminus V \right)} {\Cpc_{B(x,2r)} \left( \overline{B(x,r)}  \right)}\ge \eta \hspace{2mm}\forall x \in V}.
\ee
}\end{definition}

Note that $w_\eta(V)$ is an increasing function of $\eta \in (0,1)$ and is also an increasing 
function of the set $V$.

\begin{lemma}\label{l:cw} (See \cite[  (2.1)]{Aik01} and
\cite[Lemma 4.12]{GyS}) There exists 
$\eta=\eta(c_U, C_U, C_H)  \in (0,1)$ and $A_1>0$ such that 
	\[
	w_\eta\left( \Sett{x \in U}{ \delta_U(x)   < r} \right) \le A_1 r.
	\]
\end{lemma}

\proof Set $V_r=\Sett{x \in U}{ \delta_U(x)  < r}$. 
By Lemma \ref{l:faraway}
there is a constant $A_1>1$ such that for any point 
$x \in V_r$, there is a point $z \in U \cap B(x,A_1 r)$ with the property that 
$\delta_U(z) > 2r$. 
 By domain monotonicity of capacity, we have
	\[
	\Cpc_{B(x,2A_1r)} \left( \overline{B(x,A_1r)} \setminus V_r \right) \ge \Cpc_{B(x,2A_1r)} \left( \overline{B(z,r)} \right) \ge \Cpc_{B(z,3A_1r)} \left( \overline{B(z,r)} \right). 
	\]
 The capacities $\Cpc_{B(z,3A_1r)} \left( \overline{B(z,r)} \right)$ and $\Cpc_{B(x,2A_1r)} \left( \overline{B(x,A_1r)} \right) $ are comparable by Proposition \ref{p-consq}(a)-(c), and so the condition in
 \eqref{e:w-eta} holds for some $\eta>0$, with $r$ replaced by $A_1r$.
  \qed

We now fix $\eta \in (0,1)$ once and for all, small enough such that the 
conclusion of Lemma \ref{l:cw} applies. 
In what follows, we write $f \asymp g$, if there exists a constant 
$C_1= C_1(c_U, C_U, C_H)$ 
such that $C_1^{-1}g \le f \le C_1 g$. 


\begin{definition}
{\rm Let $(X_t, t\ge 0, \bP^x, x \in \sX)$ be the Hunt process
associated with the MMD space $(\sX,d,\mu,\sE,\sF)$.
For a Borel subset $U \subset \sX$ set 
\be  
\tau_U :=  \inf \Sett{t >0} { X_t \notin U}.
\ee
Let $\Omega \subset \sX$ be open and relatively compact in $\sX$.
Since the process $(X_t)$ is continuous,
$X_{\tau_\Omega}\in \partial\Omega$ a.s.
We define the \emph{harmonic measure} $\omega(x,\cdot,\Omega)$
on $ \partial \Omega$   
by setting
\[\omega(x,F,\Omega):=   \bP^x( X_{\tau_\Omega}\in F)
\hbox{ for } F \subset  \partial \Omega. 
\]
}\end{definition}
The following lemma provides an useful estimate of the harmonic measure in terms of the capacitary width.
\begin{lemma}\label{l:hm1}(See \cite[Lemma 1]{Aik01}, and \cite[Lemma 4.13]{GyS}) There exists $a_1 \in (0,1)$ such that
for any non-empty open set $V \subset \sX$ and for all $x \in \sX$, $r>0$, 
\[
\omega\left(x, V \cap \partial B(x,r), V \cap B(x,r) \right) \le \exp\left( 2 - \frac{a_1 r}{w_\eta(V)} \right).
\]
\end{lemma}

\proof
	The proof is same as \cite[Lemma 4.13]{GyS} except that
we use Proposition \ref{p-consq}(a),(b) instead of \cite[Lemma 4.8]{GyS}. \qed
 
In the following lemma, we provide an upper bound of the harmonic measure 
in terms of the Green function. It is an analogue of \cite[Lemma 2]{Aik01}.

\begin{lemma}\label{l:hm2}(Cf. \cite[Lemma 4.14]{GyS}, \cite[Lemma 4.9]{LS} and \cite[Lemma 5.3]{Lie})
There exists $A_2, C_4 \in (0,\infty)$ such that for all $r>0$, 
$\xi \in \partial_{\widetilde{U}}U$, there  exist  $\xi_{r},\xi_{r}' \in U$ that satisfy  
$d_U(\xi,\xi_{r})=4r$, $\delta_U(\xi_{r}) \ge 2 c_U r$, $d(\xi_{r},\xi_{r}')= c_U r$ and
\[
  \omega\left(x, U \cap \partial_{\wt U} B_U(\xi,2r),B_U(\xi,2r)\right) \le 
	C_4 \frac{g_{B_U(\xi, A_2 r)}(x,\xi_{r})}{g_{B_U(\xi, A_2 r)}(\xi'_{r},\xi_{r})}, 
	\q  \forall x \in B_U(\xi,r).
\]
\end{lemma}

\proof
Let  $\xi \in \partial_{\widetilde{U}}U$ and $r>0$.
Fix $A_2 \ge 2 (12+ C_U)$ so that all $(c_U,C_U)$-inner uniform curves that connect 
two points in $B_U(\xi,12r)$ stay inside $B_U(\xi,A_2 r/2)$. 
Fix  $\xi_r, \xi'_r \in U$ satisfying the given hypothesis: these points exist by Lemma \ref{l:faraway}.
Define 
\[
g'(z)=g_{B_U(\xi,A_2 r)}(z,\xi_{r}), \, \, \mbox { for $z \in B_U(\xi,A_2 r)$ }.
\]
Set $s=\min(c_U r,5r/C_U)$. Note that $B_U(\xi_r,s) = B(\xi_r, s) \subset U$.
Since $B(\xi_{r}, s) \subset B_U(\xi, A_2 r) \setminus  B_U(\xi, 2r)$, 
 using the maximum principle given by Assumption \ref{a-green}(iv) we have
\[
g'(y) \le \sup_{z \in \partial B(\xi_{r},s)} g'(z)
\text { for all } y \in B_U(\xi,2r).
\]
By Proposition \ref{p-consq}(a), we have
\[
\sup_{z \in \partial B(\xi_{r},s)} g'(z) \asymp g'(\xi'_{r}),
\]
and hence there exists $\eps_1>0$ such that 
\[
\eps_1 \frac{g'(y)}{g'(\xi'_{r})} \le \exp(-1) \q \forall y \in B_U(\xi,2r).
\]
For all non-negative integers $j$, define
\[
U_j:= \Sett{x \in B_U(\xi,A_2 r)}{\exp\left( - 2^{j+1} \right) \le  \eps_1 \frac{g'(x)}{g'(\xi'_{r})} 
< \exp\left( - 2^{j} \right)},
\]
so that $B_{U}(\xi,2r)= \bigcup_{j \ge 0} U_j \cap B_U(\xi,2r)$.
Set $V_j= \bigcup_{k \ge j} U_k$. 
We claim that there exist $c_1, \sigma \in (0,\infty)$ such that for all $j \ge 0$
\begin{equation}\label{e:te1}
w_\eta \left( V_j \cap B_U(\xi,2r) \right) \le c_1 r \exp \left( -2^j/\sigma \right).
\end{equation}

Let $x$ be an arbitrary point in $V_j \cap B_U(\xi,2r)$. 
Let
$z$ be the first point in the inner uniform curve from $x$ to $\xi_r$ which is on $\partial_U B_U(\xi_r, c_U r)$. Then by Lemma \ref{l:qhmetric}  there exists a Harnack chain of balls in 
$B_U(\xi,A_2r) \setminus \set{\xi_r}$ connecting $x$ to $z$ of length at 
most $c_2 \log\left(1 + c_3 r/ \delta_U(x) \right)$ for some constants 
$c_2,c_3 \in (0,\infty)$. 
Hence, there are constants $\eps_2, \eps_3, \sigma>0$ such that
\[
\exp (-2^j)   > \eps_1 \frac{g'(x)}{g'(\xi_r')} \ge \eps_2 \frac{g'(x)}{g'(z)}  
\ge \eps_3 \left( \frac{ \delta_U(x) }{ r } \right)^{\sigma}.
\]
The first inequality above follows from definition of $V_j$, the second follows from Proposition \ref{p-consq}(a) and the last one follows from Harnack chaining.
Therefore, we have
\[
V_j \cap B_U(\xi,2r) \subset \Sett{ x\in U}{ \delta_U(x) \le \eps_3^{-1/\sigma} r\exp\left( -2^j / \sigma\right) },
\]
which by Lemma \ref{l:cw} immediately implies \eqref{e:te1}.

Set $R_0 = 2r$ 
and for $j \geq 1$,
\[ R_j = \left( 2 - \frac{6}{\pi^2} \sum_{k=1}^j \frac{1}{k^2} \right) r. \]
Then $R_j \downarrow r$ and as in \cite{GyS}
\begin{align} \label{e:te2}
\sum_{j=1}^{\infty}  \exp \left( 2^{j+1} - \frac{ a_1 (R_{j-1} - R_j) }{ c_2 r \exp(-2^j / \sigma) } \right)
& < C < \infty;
\end{align}
here $C$ depends only on $\sigma$, $c_2$ and the constant $a_1$ in Lemma \ref{l:hm1}.

Let $\omega_0(\cdot)= \omega\left(\cdot, U \cap  \partial_{\wt U} B_U(\xi,2r), B_U(\xi,2r) \right)$ and 
set
\[
 d_j = \begin{cases} 
 	\sup_{ x \in U_j \cap B_{{U}}(\xi,R_j)} \frac{g'(\xi_r') \omega_0(x) }{ g'(x)}, \quad 
 	& \textrm{if } U_j \cap B_{{U}}(\xi, R_j) \neq \emptyset, \\
 	0,  & \textrm{if } U_j \cap B_{{U}}(\xi, R_j) = \emptyset.
 \end{cases}\]
 It suffices to show that $ \sup_{j \ge 0} d_j \le C_1< \infty$, and 
this is 
proved by iteration exactly as in \cite[Lemma 4.9]{LS} or \cite[Lemma 5.3]{Lie}. 
The only difference is that we replace $r^2/V(\xi,r)$ in \cite{LS} (or  $\Psi(r)/V(\xi,r)$ in \cite{Lie}) by $g'(\xi'_r)$. 
%
\qed 
\sms

\def\pwU{\partial_{\widetilde{U}}}

By using a balayage formula (cf. \cite[Proposition 4.3]{Lie}) and a standard argument 
(cf. \cite[pp. 75-76]{GyS}, \cite[Theorem 5.2]{Lie}), the proof of Theorem \ref{T:BHP} 
 reduces to the following estimate on the Green function. 

\begin{theorem} \label{T:GR}
(See \cite[Lemma 3]{Aik01}.) 
There exist $C_1, A_4, A_3 \in (1,\infty)$ such that for 
all $\xi \in \partial_{\wt U} U$ and for all $r >0$, we have,
writing $D= B_U(\xi,A_4r)$,
\[
\frac{g_D(x_1,y_1)}{g_D(x_2,y_1)} \le C_1 \frac{g_D(x_1,y_2)}{g_D(x_2,y_2)}
\text{ for all $x_1,x_2 \in B_U(\xi,r)$, $y_1,y_2 \in  U \cap \pwU B_U(\xi,A_3r)$.}
\]
\end{theorem}

Our proof follows Aikawa's approach, 
replacing the use of  bounds on the  Green function with Proposition \ref{p-consq}. 
However, as we are working on domains in a metric space 
rather than $\bR^d$, we need to be careful with Harnack chaining. 
On a general metric space one cannot control the length of a Harnack 
chain in a punctured domain $D  \setminus \{z\}$ by the length of a Harnack 
chain in $D$, as is done in \cite[(2.15)]{ALM}. 
For a general inner uniform domain $D$ on a metric space, 
the domain $D \setminus \{z\}$ need not even be connected.
Since this is the key argument in this paper, we provide the full proof, and
as the proof is long we split it into several Lemmas.

We define
\be \label{e:A8def}
  A_3 = \max( 2 + 2 c_U^{-1}, 7). 
\ee

\begin{lemma} \label{L:hchain}
Let $\xi \in  \partial_{\widetilde{U}}U$, $r>0$, 
and $y_1, y_2 \in  U \cap \pwU B_U(\xi,A_3r)$. 
If $\gam$ is  a $(c_U, C_U)$-inner uniform curve from $y_1$ to $y_2$ in $U$, then
$\gam \cap B_U(\xi,2 r)=\emptyset$ and $\gamma \subset \overline{B_U\left(\xi, A_3(C_U+1) r \right)}$.
\end{lemma}
 
\proof Let $z \in \gam$. If $d_U(y_1,z) \wedge d_U(y_2,z) \le (A_3-2) r$,
then by the triangle inequality $d_U(z ,\xi) \ge 2r$. If 
$d_U(y_1,z) \wedge d_U(y_2,z) > (A_3-2) r$, then using the
inner uniformity of $\gam$,
$$ \delta_U(z) \ge c_U \left(d_U(y_1,z) \wedge d_U(y_2,z) \right) >
c_U(A_3-2) r \ge 2r,  $$
which implies that $z \not\in B_U(\xi,2 r)$.  

For the second conclusion, note that for all $z \in \gamma$, 
\[
d_U(\xi,z) \le  A_3 r + \min\left( d_U(y_1,z), d_U(y_2,z) \right) \le A_3 r + L(\gamma)/2 \le A_3(C_U+1) r.
\]
\qed

For $\xi \in  \partial_{\widetilde{U}}U$ choose
$x^*_\xi \in U \cap \pwU B_U(\xi, r)$ and $y^*_\xi \in U \cap \pwU B_U(\xi,A_3 r)$
such that $\delta_U( x^*_\xi) \ge c_U r$ and $\delta_U( y^*_\xi) \ge  A_3 c_U r$.
Note that we have
\be \label{e:ystar}
\delta_U(y^*_\xi)\ge A_3 c_U r > 2r, 
\ee
so that $B(y^*_\xi,2r) \subset U$.
Let $\gam_\xi$ be an inner uniform curve from $y^*_\xi$ to $x^*_\xi$, 
and let $z^*_\xi$ be the last point of this curve which is on  $\pd B(y^*_\xi, c_Ur)$. 
We will write these points as $x^*, y^*, z^*$ when the choice of the boundary point $\xi$  is clear.  

Define
\[
A_4= A_2+C_U(A_3+  \frac{1}{4}c_U^2  +8 ).  
\]
To prove Theorem \ref{T:GR} it is sufficient to prove that, writing
$D= B_U(\xi, A_4 r)$, we have
for all $x \in B_U(\xi,r)$ and for all $y \in U \cap \pwU B_U(\xi,A_3 r)$ 
\begin{equation} \label{e:gc1}
g_D(x,y) \asymp \frac{g_D(x^*,y)}{g_D(x^*,y^*)} g_D(x,y^*).
\end{equation}

\begin{lemma} \label{L:far-y}
Let $\xi \in \partial_{\widetilde{U}}U$, $r>0$ and
let $D= B_U(\xi, A_4 r)$.
If $x \in B_U(\xi,r)$ and $y \in U \cap \pwU B_U(\xi,A_3 r)$ 
with $\delta_U(y) \ge \fract14 c_U^2 r$, then \eqref{e:gc1} holds.
\end{lemma}

\proof
Fix $x \in B_U(\xi,r)$. 
Set
$$ u_1(y') = g_D(x,y') , \q v_1(y')=  \frac{g_D(x^*,y')}{g_D(x^*,y^*)} g_D(x,y^*). $$
The functions $u_1$ and $v_1$
are harmonic in $D  \setminus \{x,x^*\}$, vanish quasi-everywhere on the boundary of $D$,
and satisfy $u_1(y^*)=v_1(y^*)$.
Let $\gam$ be a $(c_U, C_U)$-inner uniform curve from $y$ to $y^*$; by Lemma \ref{L:hchain} 
this curve is contained in $U\setminus B_U(\xi,2r)$.
So by Lemma \ref{L:goodcurve} 
$\delta_U(z) \ge \half  c_U \left(\delta_U(y) \wedge \delta_U(y^*) \right) \ge  \fract18 c_U^3 r$
for $z \in \gamma$. 
Thus we can find a Harnack chain of balls in $U \setminus \{x,x^*\}$
of radius $ \fract18 c_U^3 r$ with length less than $C=C(c_U, C_U, C_H)$
which connects $y$ and $y^*$. Therefore,  \eqref{e:gc1}  follows
from \eqref{e:hchain}. \qed


\begin{lemma} \label{L:glb}
Let $\xi \in \partial_{\widetilde{U}}U$, $r>0$, and let $D= B_U(\xi, A_4 r)$.
If $x \in B_U(\xi,r)$ and 
$y \in U \cap \pwU B_U(\xi,A_3 r)$  with $\delta_U(y) < \fract14 c_U^2 r$, then
\begin{equation} \label{e:gc1a-lb}
g_D(x,y) \ge c  \frac{g_D(x^*,y)}{g_D(x^*,y^*)} g_D(x,y^*).
\end{equation}
\end{lemma}

\proof Fix $y$ 
and call $u$ (respectively, $v$) the left-hand 
(resp. right-hand) side of \eqref{e:gc1a-lb}, viewed as a function of $x$. 
By Assumption \ref{a-green},
$u$ is harmonic in $D \setminus \set{y}$ and $v$ is harmonic in $D \setminus \set{y^*}$. 
Moreover, both $u$ and $v$ vanish quasi-everywhere on the boundary of $D$, and $u(x^*)=v(x^*)$.

Let $\gam_\xi$ and $z^*$ be as defined above.
By Lemma \ref{L:goodcurve}, we have
$\delta_U(z) \ge \half c_U \delta_U(x^*) \ge \half c_U^2 r$ for all $z \in \gam_\xi$,
and so this curve lies a distance at least $\fract14 c_U^2 r$ from $y$.
By the choice of $z^*$ the part of the curve from $z^*$ to $x^*$ lies outside
$B_U(y^*,\half c_Ur)$. Thus there exists $N_1=N_1(c_U, C_U)$ such that there is a
$2$-Harnack chain of balls in $U \setminus \{y, y^*\}$ connecting $x^*$ with $z^*$
of length at most $N_1$.
Using this we deduce that there exists $C<\infty$ such that
\be \label{e:511a}
 C^{-1} v(z^*) \le v(x^*) \le C v(z^*), \q
  C^{-1} u(z^*) \le u(x^*) \le C u(z^*).  
\ee 
Since $B(y^*,2c_U r) \subset U  \setminus \{y\}$, we can use  the EHI
and  Proposition \ref{p-consq} to deduce that 
$$  C^{-1} v(z^*) \le v(z) \le C v(z^*), \q
  C^{-1} u(z^*) \le u(z) \le C u(z^*) \text { for all } z \in  \pd B_U(y^*,c_U r). $$
Thus there exist $c_1, c_2$ such that 
\be \label{e:511b}
c_1  u(z)  \ge u(x^*) = v(x^*) \ge  c_1 c_2 v(z) \text{ for all } z \in  \pd B_U(y^*,c_U r).
\ee 
Using Lemma \ref{L:maxp} it follows that $u \ge c_2 v$ on
$U \setminus B_U(y^*,c_U r)$, proving  \eqref{e:gc1a-lb}. \qed

For $\xi \in \pd_{\wt U} U$ set
$$ F(\xi) = B_U( \xi, (A_3+3) r ) \setminus B_U( \xi, (A_3-3) r ) . $$
Let  $A_5 = A_3 + A_4. $

\begin{lemma} \label{L:zF1}
Let $\xi \in   \partial_{\widetilde{U}}U$,
and $D = B_U(\xi, A_4 r)$. 
Then
\begin{equation}\label{e:gc6}
g_D(x,z) \le C_1 g_D(x,y^*), \q \text{for all } 
x \in B_U(\xi,2r), \, z \in  F(\xi). 
\end{equation}
\end{lemma} 

\proof 
We begin by proving that
\begin{equation}\label{e:gc6a}
g_D(x,y) \le C_1 g_D(x^*,y^*), \q \text{for all } 
x \in B_U(\xi,2r), \, y \in  F(\xi). 
\end{equation}
Let $x \in B_U(\xi,2r)$, $y \in F(\xi)$.  Let 
$\widetilde{C_U}$ be the constant from Lemma \ref{l:ballcomp}. 
We have
$D \subset B(y^*, A_5 r)$, and therefore by domain montonicity of the Green 
function and Proposition \ref{p-consq}
we have for any $z\in D$ with $d(x,z) \ge r/(2 \widetilde{C_U})$,
\be \label{e:gc7}
 g_D(x,z) \le g_{B(y^*,A_5r)}(x,z) \le  g_{B(y^*,A_5r)}(x^*,y^*)
\le C_1 g_D(x^*,y^*). 
\ee
If $d(x,y) \ge r/(2 \widetilde{C_U})$ this gives \eqref{e:gc6a}.
%

\begin{figure}[h] 
	\centering
\includegraphics[width=0.5\textwidth]{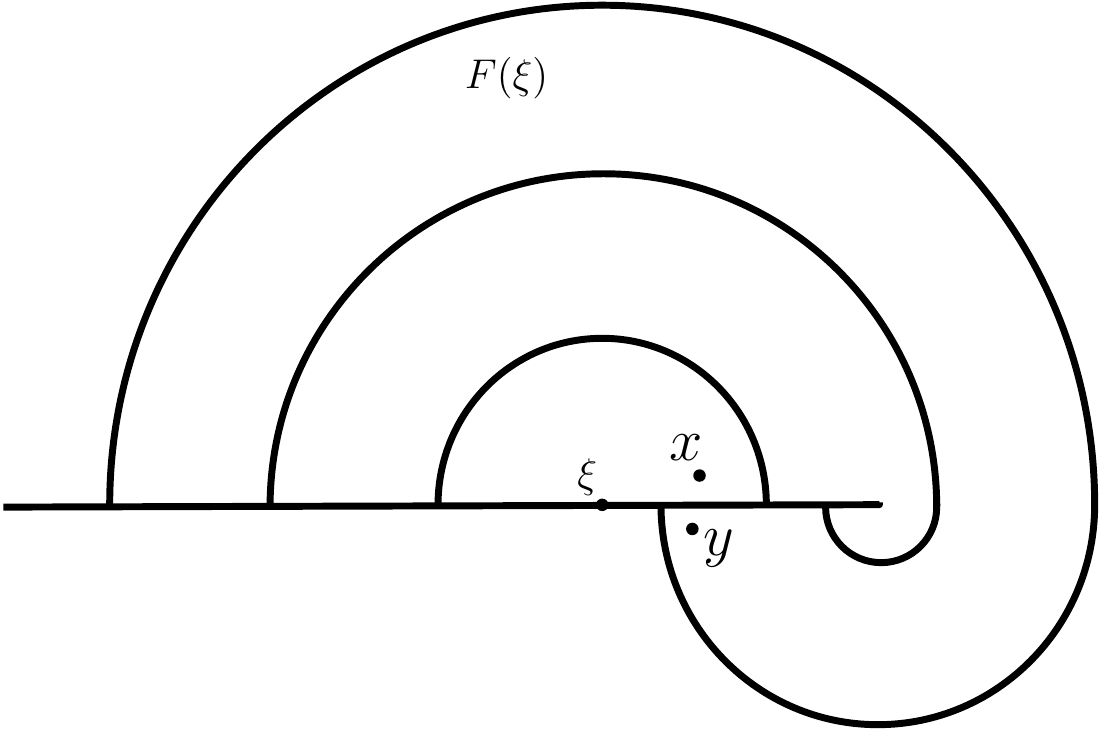}
\caption{The inner uniform domain $U= \bR^2 \setminus \left([-1,0] \times \set{0} \right)$ showing the set $F(\xi)$}
\label{f:slit}
\end{figure}

Next, we consider the case $d(x,y)< r/(2 \widetilde{C_U})$. 
(See Figure \ref{f:slit} for an example of a slit domain containing such points). 
Let $B_y$ denote the connected component of 
$p^{-1} \left(B(p(y),r/\widetilde{C_U}) \cap \overline{U}\right)$ that contains $y$. 
By Lemma \ref{l:ballcomp}, we have $B_y \subset B_{\wt{U}}(y,r)$. 
As $g_D(x,\cdot)$ is harmonic in $B_y \cap U$, by the maximum principle,  we have
\[
g_D(x,y) \le \sup_{z \in U \cap \partial_{\wt{U}}B_y}  g_{D}(x,z) \le \sup_{z \in \partial B(y,r/\widetilde{C_U})}  g_D(x,z).
\]
By the triangle inequality, 
we have $d(x,z) \ge r/(2 \widetilde{C_U})$ for all $z \in \partial B(y,r/\widetilde{C_U})$,
and therefore  \eqref{e:gc6a} follows from \eqref{e:gc7}.
This completes the proof of \eqref{e:gc6a}. 

By the continuity of the Green function,  we can extend \eqref{e:gc6a} as follows:
\begin{equation}\label{e:gc6as}
g_D(x,y) \le C_1 g_D(x^*,y^*), \q \text{for all } 
x \in U \cap \overline{B_U(\xi,2r)}^{d_U}, \, y \in  F(\xi). 
\end{equation}
Now, let $x \in B_U(\xi,2r)$, $z \in  F(\xi)$.
Since $g_D(\cdot,z)$ is harmonic in $D \setminus \set{z}$, by 
the maximum principle we have
\begin{equation} \label{e:mark1}
g_D(x,z) \le  \omega(x, U \cap \partial_{\wt U} B_U(\xi,2r),  B_U(\xi,2r))
\sup_{ x' \in U \cap \pd_{\wt{U}} B_U(\xi,2r)} g_D(x',z) . 
\end{equation}
We use Lemma \ref{l:hm2} to bound the first term, and \eqref{e:gc6as}
to bound the second, and obtain
\be \label{e:gr1}
g_D(x,z) 
\le c \frac{g_{B_U(\xi, A_2 r)}(x,\xi_{r})}{g_{B_U(\xi, A_2 r)}(\xi'_{r},\xi_{r})}
 g_D(x^*,y^*).
\ee

We then have by Proposition \ref{p-consq}(a)-(c),  Harnack chaining, and domain monotonicity
$$  g_{B_U(\xi, A_2 r)}(\xi'_{r},\xi_{r}) \asymp g_D(x^*,y^*),
\quad 
 g_{B_U(\xi, A_2 r)}(x,\xi_{r}) \le c g_D(x,y^*), $$
and combining these inequalities completes the proof of \eqref{e:gc6}. Note that for the second inequality above, one needs to consider two different cases: $\delta_U(x) \le \half c_U^2 r$ and  $\delta_U(x) > \half c_U^2 r$. 
\qed

\begin{lemma} \label{L:gub}
Let $\xi \in \partial_{\widetilde{U}}U$,
and $D = B_U(\xi, A_4 r)$.
If $x \in B_U(\xi,r)$ and 
$y \in \partial_U B_U(\xi,A_3 r)$  with $\delta_U(y) < \fract14 c_U^2 r$, then
\begin{equation} \label{e:gc1a-ub}
g_D(x,y) \le c  \frac{g_D(x^*,y)}{g_D(x^*,y^*)} g_D(x,y^*).
\end{equation}
\end{lemma}

\proof
Let $\zeta \in \partial_{\widetilde{U}} U$ be a point such that $d_U(y,\zeta) < c_U^2 r/4$,
and let $\zeta_r$ and $\zeta'_r$ be the points given by Lemma \ref{l:hm2}
corresponding to the boundary point $\zeta$.  
Since $g_D(x,\cdot)$ is harmonic in $B_{U}(\zeta,2r)$, we have
\begin{equation} \label{e:mark2}
g_D(x,y) \le  \omega(y, \partial_U B_U(\zeta,2r),  B_U(\zeta,2r))
\sup_{ z \in U \cap  \partial_{\wt U} B_U(\zeta,2r)} g_D(x,z) . 
\end{equation}
Since $B_U(\zeta,2r) \subset F(\xi)$,  by Lemma \ref{L:zF1}, 
the second term in \eqref{e:mark2} is bounded by $c g_D(x,y^*)$.
Using Lemma \ref{l:hm2} to control the first term, we obtain
\be
g_D(x,y) \le  c g_D(x,y^*) \frac{g_{B_U(\zeta,A_2 r)}(y, \zeta_r) }
   {g_{B_U(\zeta,A_2 r)}(\zeta'_r, \zeta_r) }.
\ee
Again by Harnack chaining, Proposition \ref{p-consq}, and domain monotonicity we have
$$   g_{B_U(\zeta,A_2 r)}(\zeta'_r, \zeta_r) \asymp g_D(x^*,y^*), $$
and
$$ g_{B_U(\zeta,A_2 r)}(y, \zeta_r) \le c g_D(y, x^*), $$ 
and combining these estimates completes the proof. \qed

\sm 
{\em Proof of Theorem \ref{T:GR}. }
The estimate \eqref{e:gc1} follows immediately from Lemmas \ref{L:far-y}, \ref{L:glb}
and \ref{L:gub}, and as remarked before, the Theorem follows 
from \eqref{e:gc1}. \qed

\begin{remark} 
{\rm One might ask if the converse to Theorem \ref{T:BHP} holds. 
That is, suppose $(\sX,d,\mu,\sE,\sF)$ is a MMD space such that for every
inner uniform domain the BHP holds. Then does the EHI hold for
 $(\sX,d,\mu,\sE,\sF)$?
 
The following example shows this is not the case.
Consider the measures  $\mu_\alpha$ on $\bR$ 
 given by $\mu_\alpha (dx) = (1+ \abs{x}^2)^{\alpha/2} \, \lambda(dx)$,  where $\lambda$ denotes the Lebesgue measure. 
(See  \cite{GS}). The Dirichlet forms 
\[
\sE_\alpha(f,f) =  \int_{\bR} \abs{ f'(x)}^2 \mu_\alpha(dx)
\]
do not satisfy the Liouville property if $\alpha >1$. 
This is because the two ends at $\pm \infty$ are transient,
so the probability that the diffusion eventually ends up in $(0,\infty)$ is a non-constant positive harmonic function. Since the Liouville property fails, so does the EHI.

On the other hand, the space of inner uniform domains in $\bR$ is same 
as the space of (proper) intervals in $\bR$.
The space of harmonic functions in a bounded interval vanishing at a boundary point is one dimensional, and hence the BHP holds. We can take $R(U)$ in Theorem \ref{T:BHP} 
as $\diam(U)/4$. In view of this example, the following question remains open: \emph{Which diffusions admit the scale invariant BHP for all inner uniform domains?} Theorem \ref{T:BHP} shows 
that the EHI provides a sufficient condition for the scale invariant BHP,
but the example above shows that the EHI is not necessary. 
} \end{remark}

We now give two examples to which Theorem \ref{T:BHP} applies 
but earlier results do not.

\begin{example} \label{ex:new}
{\rm (1) (See \cite[Example 6.14]{GS}) 
Let $n \ge 2$. Consider the measure 
$\mu_\alpha(dx)= \left(1 + \abs{x}^2\right)^{-\alpha/2}\,\lambda(dx)$, 
where $\lambda$ is  Lebesgue measure on $\bR^n$. 
The second order `weighted Laplace' operator $L_\alpha$ on $\bR^n$
associated with the measure $\mu_\alpha$ is given by
\[
L_\alpha= \left(1 + \abs{x}^2\right)^{-\alpha/2} \sum_{i=1}^n \frac{\partial}{\partial x_i} \left(\left(1 + \abs{x}^2\right)^{\alpha/2} \frac{\partial}{\partial x_i} \right)= \Delta  +  \alpha\frac{ x. \nabla}{1 + \abs{x}^2}.
\]
The operator $L_\alpha$ is the generator of the Dirichlet form 
\[
\sE_\alpha(f,f)= \int_{\bR^n} \norm{\nabla f}^2 \, d\mu_\alpha,
\]
on $L^2(\mathbb{R}^n, \mu_\alpha)$.
Grigor'yan and Saloff-Coste \cite{GS} show that $L_\alpha$ satisfies the PHI if and only if $\alpha> -n$ but satisfies 
the EHI for all $\alpha \in \bR$. If $\alpha \le n$, the measure $\mu_\alpha$ does \emph{not} satisfy the volume doubling property.
Assumption \ref{a-green} for this example follows from Lemmas \ref{l-ex-philoc}(a) and \ref{l-ex-adir}.

\noindent (2) The first example of a space that satisfies the EHI but fails to satisfy the volume doubling 
property was given by Delmotte \cite{Del}, in the graph context.
A general class of examples similar to \cite{Del} is given in \cite[Lemma 5.1]{Bar}.
The associated cable systems of these graphs do satisfy the EHI, but
do not satisfy a global parabolic Harnack inequality of the kind given in
Definition \ref{d-locvdpiphi}, i.e. $(\operatorname{PHI(\Psi)})_{\operatorname{loc}}$
with $R=\infty$.

}\end{example}

\bibliographystyle{amsalpha}

\noindent MB: Department of Mathematics,
University of British Columbia,
Vancouver, BC V6T 1Z2, Canada. \\
barlow@math.ubc.ca

\noindent MM:  Department of Mathematics, University of British Columbia and Pacific institute
for the Mathematical Sciences, 
Vancouver, BC V6T 1Z2, Canada. \\
mathav@math.ubc.ca

\end{document}